\documentclass[12pt, leqno]{article}
\usepackage{amsmath}
\usepackage{amsfonts}
\usepackage{amssymb}
\usepackage{amssymb,esint}
\usepackage{graphicx,epsfig,color,colordvi}
\usepackage{color}
\usepackage{color}
\newtheorem{thm}{Theorem}[section]
\newtheorem{cor}[thm]{Corollary}

\newtheorem{prop}[thm]{Proposition}

\numberwithin{equation}{section}

\newcommand{\RR}{\mathbb{R}}
\newcommand{\re}{\mathbb{R}}
\newcommand{\ren}{\mathbb{R}^N}

\newcommand{\ve}{\varepsilon}

\def\qed{\,\unskip\kern 6pt \penalty 500
\raise -2pt\hbox{\vrule \vbox to8pt{\hrule width 6pt
\vfill\hrule}\vrule}\par}
\definecolor{darkblue}{rgb}{0.05, .05, .65}
\definecolor{darkgreen}{rgb}{0.05, .70, .05}
\definecolor{darkred}{rgb}{0.8,0,0}
\def\qed{\unskip\kern 6pt \penalty 500
\raise -2pt\hbox{\vrule \vbox to8pt{\hrule width 6pt
\vfill\hrule}\vrule}\par}
\begin{document}
\title{\textbf{Barenblatt  solutions and asymptotic \\ behaviour for a  nonlinear fractional heat\\ equation of porous medium type }}
\author{\Large Juan Luis V\'azquez \footnote{ Dpto. de Matem\'aticas, Univ. Aut\'onoma de Madrid, 28049 Madrid, Spain. \newline   E-mail: {\tt juanluis.vazquez@uam.es}}
 \\
{ \ }
\\\rightline{\small \sl This paper is dedicated to Grisha Barenblatt,}
\\ \rightline{\small \sl maestro and friend, for his 85th birthday}
}

\date{ }

\maketitle

\begin{abstract}
We establish the existence, uniqueness and main properties of the fundamental solutions for the fractional porous medium equation introduced in \cite{PQRV1}. They are self-similar functions of the form \ $u(x,t)= t^{-\alpha} f(|x|\,t^{-\beta})$ \
with suitable $\alpha$ and $\beta$. As a main application of this construction, we prove that the asymptotic behaviour of general solutions is represented by such special solutions. Very singular solutions are also constructed. Among other interesting qualitative properties of the equation we prove an Aleksandrov reflection principle.
\end{abstract}

\

\vskip 1cm


\medskip


\tableofcontents

\newpage
\section{Introduction}
\label{sec.intro}

Our main goal is to determine the existence, uniqueness and main properties of the solution  of the equation
\begin{equation}  \label{eq1}
\partial_tu+(-\Delta)^{s}(u^m)=0,  \qquad  x\in\mathbb{R}^N,\; t>0\,,
\end{equation}
with fractional exponent $s\in (0,1)$, and  $m>0$. This equation, or better said family of equations, is one of the standard models of nonlinear diffusion involving long-distance effects in the form of fractional Laplacian operators, which are the most representative nonlocal operators of elliptic type. We recall that the fractional Laplacian operator is a kind of differentiation operator of order $2s$, for arbitrary $s\in (0,1)$, that can be conveniently defined  through its Fourier Transform symbol, which is $|\xi|^{2s}$, so that for $s=1$ we recover the standard Laplacian. A major difference between the standard and the fractional Laplacian is best seen in the stochastic point of view, and it consists  in taking into account long-range interactions in the latter, which explains features that we will see below, like enhanced propagation with the appearance of fat tails at long distances (to be compared with the typical exponentially small tails of the standard diffusion or the compactly supported solutions of porous medium flows). This is known as anomalous diffusion.  There is a wide literature on the subject, both for its relevance to Analysis, PDEs, Potential Theory, Stochastic Processes, and for the growing number of applications in Mechanics and other applied fields.  See basic information in \cite{AbeThurner05, CS07, Landkof, Stein70, Valdinoci,VazAbel}.

We consider the Cauchy problem for Equation \eqref{eq1} taking as initial data a Dirac delta,
\begin{equation}\label{eq.id}
u(x,0) = M\delta(x)\, \qquad M>0.
\end{equation}
Solutions with such data are called fundamental solutions in the linear theory, and we will keep that name though their relevance is different in the nonlinear context.
 We use the concept of  continuous and nonnegative weak solution introduced in \cite{PQRV1}, \cite{PQRV2}, for which there is a well-developed theory when the datum is a  function in $L^1(\ren)$.

Putting $s=1$ we recover the standard Porous Medium Equation (PME), where the question under discussion is well-known, see the comments in the next section. The exponent $m$ varies in principle in the range $m>1$, but the methods extend to the linear case $m=1$. and even to the fast diffusion range (FDE)  $m<1$ on the condition that $m> m_c=\max\{(N-2s)/N,0\}$. Such type of restriction on $m$ from below carries over from the PME-FDE theory, \cite{Vazsmooth06}.

A main result of this paper consists in proving that there exists a  unique fundamental solution of Problem \eqref{eq1}--\eqref{eq.id} as follows.

\begin{thm} \label{thm.Bs} For every choice of parameters $s\in(0,1)$ and $m>m_c=\max\{(N-2s)/N,0\}$, and every $M>0$, equation \eqref{eq1} admits a unique fundamental solution; it is a nonnegative and continuous weak solution for \ $t>0$ and takes the initial data \eqref{eq.id} as a trace in the sense of Radon measures. Such solution has the self-similar form
\begin{equation}\label{sss}
u_M^*(x,t)=t^{-\alpha} f(|x|\,t^{-\beta})
 \end{equation}
for suitable $\alpha$ and $\beta$ that can be calculated in terms of $N$ and $s$ in a dimensional way, precisely
\begin{equation}\label{scale.expo}
\alpha=\frac{N}{N(m-1)+2s}, \qquad \beta=\frac{1}{N(m-1)+2s}\,.
\end{equation}
The profile function  $f(r)$, $r\ge 0$, is a  bounded and H\"older continuous function,\footnote{$f$ is actually a smooth function, but such result is not discussed or used here.} it is positive everywhere, it is monotone and  goes to zero at infinity.
\end{thm}


More precisely, the initial data are taken in the weak sense of measures
\begin{equation}
\lim_{t\to 0}\int_{\ren} u(x,t)\phi(x)\,dx= M\phi(0)
\end{equation}
for all $\phi\in C_b(\ren)$, the space of continuous and bounded functions in $\ren$.
We will call these self-similar solutions  of Problem \eqref{eq1}-\eqref{eq.id} with given $M>0$ the {\sl Barenblatt solutions} of the fractional diffusion model by analogy with the PME and other prominent studied cases. The form of the exponents explains the already mentioned restriction on $m$ from below. The result is proved in Sections \ref{sec.prelim} to \ref{sec.1d}.

We then prove that the asymptotic behaviour as $t\to\infty$ of the class of nonnegative
weak solutions of \eqref{eq1} with finite mass  (i.\,e., $\int u(x,t)\,dx<\infty$) is given in first approximation by the family $u_M^*(x,t)$, i.\,e., the Barenblatt solutions are the attractors in that class of solutions. See Section \ref{sec.ab},  Theorem \ref{thm.exlimit}.

In the limit $M\to\infty$ we obtain a special solution with a fixed isolated singularity at $x=0$ that has separate-variables form and we call very singular solution (VSS) by analogy  with the standard Fast Diffusion Equation. This happens precisely in the range $m_c<m<m_1=N/(N+2s)$, see Theorems \ref{thm.vss1}, \ref{thm.vss2}. The result represents a marked difference with the case $s=1$ where VSS exist in the larger range  $m_c<m<1$. It is another manifestation of the long-range interactions of the fractional Laplacian, that avoids some of the purely local estimates of the standard FDE with the classical Laplacian operator; indeed, extrapolation of the standard estimates would justify the existence of a VSS in cases there is none. Let us recall that the class of VSS has played a role in the study of nonlinear parabolic equations since the seminal paper of Brezis, Peletier and Terman \cite{BPT86}.

We then devote one Section to constructing eternal solutions in the critical exponent case $m=(N-2s)/N$. They have exponential self-similarity.

The study of evolution equations with fractional Laplacian operators implies the need to develop in this setting tools that have been successful in the standard and parabolic theory. One of them is the Aleksandrov reflection principle that we state and prove in the elliptic and parabolic settings in Section \ref{sec.alek}.

A number of appendices and a section of comments and extensions close the paper.

\section{Motivation and historical perspective}\label{sec-hist}

 The question of finding fundamental solutions of elliptic and parabolic equations is one of paramount importance in the study of linear elliptic and parabolic equations, see e.\,g. \cite{Evans, Horm}. In the nonlinear theory their role is not so apparent, but it has been shown that they can be an important tool in the existence and regularity theory, and very important in the description of the asymptotic behaviour.

 Generalizing the classical heat equation, a prominent case at hand in the evolution theories  is the Porous Medium Equation (PME), $\partial_t u -\Delta(u^m)=0$, $m>1$, introduced in the last century in connection with a number of physical applications and extensively studied as a prototype of nonlinear diffusive evolution with interesting analysis and geometry (free boundaries). Fundamental solutions were discovered in the 1950's by Zeldovich and Kompanyeets \cite{ZK50} and Barenblatt \cite{Bar52}, and later by Pattle \cite{Pattle59},  the latter made a complete description. They have been usually called Barenblatt solutions in the literature. This discovery was so to say the starting point of the rigorous mathematical theory that has been gradually developed since then.  Their role in the asymptotic behaviour of general solutions of the PME was established in papers by Kamin \cite{Kam76, Kam78},  Friedman and Kamin \cite{FrKa80} and the author \cite{VAspme03}. They have also played an important role in the existence and regularity theory in the work by many authors. The monograph \cite{vazquez07} contains  a detailed exposition of the topic.

The surprising relation between existence and uniqueness of fundamental solutions and precise asymptotic behaviour relies on the existence of a scaling group under which the solutions of the PME are invariant. It implies that a fundamental solution is in fact a self-similar function. This is what we call a {Barenblatt solution}. We will see below how the scaling group works in the case of equation \eqref{eq1}. Self-similarity plays a big role in our understanding of fundamental processes in mathematics and mechanics, as explained by Prof. Barenblatt in his books \cite{Barbk96, Barbk03}.

Following the analysis of Barenblatt solutions for the PME, other equations have been explored. To name a few, the Fast Diffusion Equation, which is $\partial_t u -\Delta(u^m)=0$, with $m<1$. Though it looks formally the same, the qualitative theory changes. The role of the Barenblatt solutions remains basically unchanged only for $m>(N-2)/N$, while for $m<(N-2)/N$ the whole functional setting changes abruptly, cf. \cite{Vazsmooth06} and \cite{BBDGV, BDGVpnas}.
Barenblatt solutions have been constructed for the $p$-Laplacian equation, $\partial_t u -\Delta_p(u)=0$ \cite{Bar52b}; in this case they exist in the range $2N/(N+1)<p<\infty$ and their role and properties are quite similar (loosely speaking) to the porous medium case. These solutions have been used extensively in the regularity theory \cite{DiBbook} and in the asymptotic behaviour \cite{KV88}. The ideas can be extended to the doubly-nonlinear diffusion equation $\partial_t u -\Delta_p(u^m)=0$, as  proposed in  \cite{Bar52b}, see \cite{ABC2010} for recent results, and to other models, not necessarily of diffusive type. Among such models we count conservation laws, \cite{LP84}.

In some cases the fundamental solutions have similarity exponents that cannot be calculated from dimensional considerations, giving rise to so-called anomalous exponents, it is also described as self-similarity of the second kind. Then the initial distribution is not a Dirac delta but a more general mass distribution, in other words some singularity located at $x=0$, $t=0$. We say that it is a {\sl solution with an isolated singularity}, or {\sl a source-type solution}, but not a fundamental solution. There is a wide literature on solutions with isolated singularities, both in elliptic and parabolic problems, among them particular attention has given to the so-called very singular solutions, cf. \cite{BenBrez, BPT86, KamPel85, KPV89, Veron}.
Again, the influence of Prof. Barenblatt has been felt in the form of interesting models from the physical sciences. One of them is the Barenblatt equation of elasto-plastic filtration, proposed in \cite{Barbk96},  where the source-type solutions have been studied in \cite{BarSiv69, KPV91}, see also \cite{HV12}. They correspond to  self-similarity of the second kind, with anomalous exponents.  Another model is the equation of turbulent bursts \cite{Bar73} treated mathematically in \cite{KV92}, \cite{HPel92} and \cite{Hulshof}, where self-similarity with anomalous exponents is also found as asymptotic behaviour of a wide class of solutions. Anomalous exponents are discussed in \cite{ArVaz95} in connection with focussing problems, and they are found in the Fast Diffusion Equation under the critical exponent \cite{Ki93, PZ95, Vazsmooth06}, but in these cases they do not correspond to source-type solutions.

\smallskip

\noindent {\sc Nonlinear diffusion with nonlocal operators.}
There are many other instances of the paradigm we are discussing. Recently there has been a remarkable interest in mathematical models of diffusion involving long-range effects represented by  fractional Laplacian  operators and other singular integral operators. This responds to a serious motivation from the physical sciences and has produced so far quite interesting mathematical developments. Let us add some comments to the information that opened the paper on fractional Laplacian operators.  For functional analysis considerations, it is usually best to consider the inverse operator $(-\Delta)^{-s}$, which happens to be the integral operator associated to the Riesz kernels, as we will do often in this paper. Also, a special interest has the case $N=1$ where the operator $(-\Delta)^{s}$ is the natural realization of the concept of $(\partial_x)^{2s}$ (up to a constant factor)  in the form of a positive symmetric operator.

Regarding nonlinear evolution models, the author has been involved in the analysis of two of such models, called fractional porous medium equations. One of them is Equation \eqref{eq1}  that we study here, while the other model uses the equation
\begin{equation}\label{eq.mod1}
    \partial_t u -\nabla \cdot (u\nabla (-\Delta)^{-s}u)=0, \quad 0<s<1\,.
\end{equation}
Both models have quite different properties, see a survey of recent results in \cite{VazAbel}.
The construction of the Barenblatt solutions for this latter model has been performed by Caffarelli and the author in \cite{CaVa11}. The asymptotic behaviour is obtained in  \cite{CaVa11}  using sophisticated entropy and obstacle problem methods. An explicit form for the Barenblatt solutions is found in \cite{bik}
\begin{equation}
U(x,t)=c_1\,t^{-\alpha} (1-c_2|x|^2 t^{-2\alpha/N})_+^{1-s}\,  \qquad \text{with} \quad \alpha=N/(N+2-2s).
\end{equation}
In the limit $s\to 1$ we obtain  quite interesting fundamental solutions of the non-diffusive limit $u_t -\nabla \cdot (u\nabla (-\Delta)^{-1}u)=0$, of interest in superconductivity and superfluidity. The solutions take the form
\begin{equation}
u_M^*(x,t)=\frac1t \chi_{B_R(t)}(x), \quad R(t)=c\,t^{1/N}\,, \quad c>0\,,
\end{equation}
as described in \cite{BLL} and \cite{SeVa12}. Each one  represents a round vortex patch that expands to fill the space as $t\to\infty$.

\medskip

\noindent {\sc The linear case.}
  Solving Equation \eqref{eq1} in the linear case $m=1$  is  easier. As explained in \cite{PQRV2}, when  $m=1$ and $0<s<1$ the evolution has  the integral representation
\begin{equation} \label{lineal}
u(x,t)=\int_{\mathbb{R}^N}K_s(x-z,t)u_0(z)\,dz,
\end{equation}
where $K_s$ has Fourier transform  $\widehat K_s(\xi,t)=e^{-|\xi|^{2s} t}$. This means that, for
$0<s<1$, the kernel $K_s$ has the form $K_s(x,t)=t^{-N/2s}F(|x|\,t^{-1/2s})$ for some profile
function $F$ that is positive and decreasing, and behaves at infinity like $F(r)\sim r^{-(N+2s)}$, see \cite{Blumenthal-Getoor}.
When $s=1/2$, the kernel is explicit,
\begin{equation}
K_{1/2}(x,t)=C_N\,t\,(x^2+t^2)^{-(N+1)/2}\,.
\end{equation}
If $s=1$ the function $K_1(x,t)$ is the Gaussian heat kernel. The linear model has been well studied by probabilists.

However, an integral representation of the evolution like \eqref{lineal} is not available in the nonlinear case; the tools to treat the nonlinear case were developed in \cite{PQRV1}, \cite{PQRV2}.
Below, we  perform the analysis of existence, uniqueness, properties and applications of Barenblatt solutions for equation \eqref{eq1}.  We will concentrate on the difficulties caused by the nonlinear form of the equation. Our Barenblatt solutions will play the role of the kernel $K_s$ in the study of asymptotic behaviour. Our results below include the linear case as a particular instance. Let us recall that the nonlinear techniques are quite different from the linear ones.

\section{Preliminaries}\label{sec.prelim}

Before we start the construction we need to review some basic facts.

\noindent {\sc Scaling.} Equation \eqref{eq1} is invariant under translations in space and time, and also under the scaling group, and this fact will play a big role in what follows. Let first observe that for every solution $u(x,t)$ and constants $A,B,C>0$ the function
\begin{equation}\label{scaling}
\widehat u(x,t)= A\, u(Bx,Ct)
\end{equation}
is again a solution of \eqref{eq1} if \ $C=A^{m-1}B^{2s}$. This generates the whole scaling group. If moreover, we impose the condition that the solutions have constant finite mass (i.\,e.,  integral in space), then $A=B^{N}$, which implies that the group is reduced to the one-parameter family
\ $\widehat u(x,t) =  B^{N} u(B x, B^{N(m-1)+2s} t)$,
which is usually written in terms of $\lambda=B^{N(m-1)+2s}$ as
\begin{equation}\label{scaling1}
(T_\lambda u)(x,t)=\lambda^{\alpha} u(\lambda^{\beta}x, \lambda\, t)\,,
\end{equation}
with scaling exponents given by formula \eqref{scale.expo}:
\begin{equation*}
\alpha=\frac{N}{N(m-1)+2s}, \qquad \beta=\frac{1}{N(m-1)+2s}\,,
\end{equation*}
Note that $\alpha$ and $\beta>0$ iff $m>m_c=(N-2s)/N$.  The values of both parameters will be fixed in what follows.

Another application of scaling consists in reducing solutions to unit mass. Indeed, if $u(x,t)$ is a solution with $\int u(x,0)\,dx=M>0$ then
\begin{equation}\label{mass.red}
\widehat u(x,t) =  M^{-1} u(x, M^{-(m-1)} t),
\end{equation}
is another solution with unit initial mass, $\int \widehat u(x,0)\,dx=1$.

\medskip

\noindent {\sc Potential equation.} Take the convolution $U(x,t) =u(x,t)\ast {\cal I}_{2s}(x)$,
 where ${\cal I}_{2s}$ is the Riesz kernel
\begin{equation}\label{riesz1}
{\cal I}_{2s}(x)= C_{N,s}|x|^{-(N-2s)}\,,
\end{equation}
that is the kernel of operator $(-\Delta)^{-s}$, $0<s<1$.\footnote{For the value of $C_{N,s}$ see formula \eqref{riesz.formula} below.}  Use of this kernel restricts the range of application to $0<s<1/2$ when $N=1$. This is a typical difficulty of potential theory that is found in the standard heat and porous medium equations in dimensions $N=1,2$. We will discuss the case $N=1$, $1/2\le s<1$, separately in Section \ref{sec.1d}.

After this caveat, we resume the theory.  We apply the kernel to a solution $u(x,t)$:
$$
U(x,t)= C_{n,s}\int \frac{u(y,t)}{|x-y|^{N-2s}}\,dy \,.
$$
Then $(-\Delta )^s U= u$, and using the equation we get equation (PE):
\begin{equation}\label{eq.pot1}
U_t=((-\Delta )^{-s} U)_t=-(-\Delta )^{-s}(-\Delta )^{s}u^m=-u^m.
\end{equation}
This is easy to justify since $u^m(t)\in L^2$ with $u_t=(-\Delta )^{s}u^m\in L^1$ for every $t>0$, according to the theory.

The scaling group acts also on the solutions of the PE: If $A,B>0$ are positive and
$U$ is a solution of the PE, then so is
\begin{equation}\label{scalPE}
\widehat U(x,t) =  A \,U(B x, Ct), \qquad \text{with } \ C=A^{m-1}B^{2m}
\end{equation}
is again a solution of the same equation with suitably rescaled initial data.

The idea of using the potential equation to prove uniqueness of solutions of diffusive equations with measure data goes back to Pierre's work for the standard PME \cite{Pierre}. This has been followed later by a number of  authors in various contexts.

\section{Existence of solutions with measure data}\label{sec.ex}

Here is the basic result in this topic. ${\cal M}^+(\ren)$ is the space of bounded and nonnegative Radon measures on $\ren$.

\begin{thm} For every $\mu\in {\cal M}^+(\ren)$ there exists a nonnegative and continuous weak solution of Equation \eqref{eq1} in $Q=\ren\times (0,\infty)$ taking initial data $\mu$ in the sense that for every $\varphi\in C_c^{2}(\ren)$ we have
\begin{equation}
\lim_{t\to 0} \int u(x,t)\varphi(x)\,dx=\int \varphi(x)d\mu(x)\,.
\end{equation}
\end{thm}

\noindent{\sl Proof.} \noindent (i) Existence comes from approximation of the initial data $\mu$  with a smooth mollifier sequence $\rho_\ve(x)$, thus we take smooth initial data $\mu_\ve=\mu\ast \rho_\ve$. Using the results of \cite{PQRV2} we know that the corresponding solutions $u_\ve$ exist and are uniformly bounded in $L^\infty(0,\infty:L^1(\ren))$ and  $L^\infty (\ren\times (\tau,\infty))$. The decay rate with time,  $ \|u(t)\|_\infty\le Ct^{-\alpha}$, is a main tool, and it holds uniformly for the whole sequence. They also have uniform bounds  on the energy $\int_t^\infty \int |(-\Delta)^{s/2} u^m|^2\,dxdt$, cf. \cite{PQRV2}, formula (8.5).

\medskip

\noindent (ii) We pass to the limit $\ve\to 0$. Convergence holds in principle up to subsequences $\ve_n\to 0$. The solutions $u_\ve$ and the  limit $\overline{u}$ are uniformly $C^\alpha$ in space and time for all $t\ge \tau>0$. The limit is easily proved to be a weak solution of the equation.

\medskip

\noindent (iii) The initial data are taken  in the sense of initial traces, since the weak formulation passes to the limit for $t\ge \tau>0$ and the second integral term is uniformly small for small $t\in (0,\tau)$: if $L=(-\Delta)^s$ then
$$
|\iint u^m L \varphi \,dxdt|\le \int_0^t \|u(s)\|_\infty^{m-1}\left(\int u\, |L\varphi|\,dx\right)\,ds\,.
$$
 Now we use that $|L\varphi|\le C$, $\int u\,dx\le C$ and the decay $ \|u(t)\|_\infty\le Ct^{-\alpha}$ with $\alpha=N/(N(m-1)+2s)$.
Therefore,
\begin{equation}
\alpha (m-1)=\frac{N (m-1)}{N(m-1)+2s}<1
\end{equation}
and the integral in time converges at $t=0$.  This calculation is valid for the $u_\ve$ and also for the limit $\overline{u}$. See a similar argument  in  \cite{SeVa12}, Section 3. So we have initial traces in the sense of distributions. Even better, since the test functions can be taken in $W^{2s,\infty}(\ren)$ for instance.

\medskip

\noindent (iv) The previous argument as $t\to 0$ works for $m\ge 1$. For $m_c<m<1$, we use good test functions and H\"older estimates:
\begin{equation}\label{m<1}
|\iint Lu^m \varphi \,dxdt|\le \left(\iint u dxdt\right)^{m}\left(\iint (L\varphi)^{1/(1-m)} \right)^{1-m}\le Ct^{1-m}.
\end{equation}
\qed

\section{Approximate data for the fundamental solution}
\label{sec.ma}

\noindent {\bf 1}. In our study of the solutions with Dirac deltas as initial data, we will use some monotone approximations. The simplest approach is to take a smooth convolution kernel $\rho\in C^\infty_c(\ren)$, $\rho\ge 0$, $\int \rho\,dx=1$. We also take $\rho$ to be radial, i.\,e., a function of $|x|$ and supported in the ball of radius 1. Put $\rho_\ve=\ve^{-N}\rho(x/\ve)$, and define the regularization $V_\ve =\rho_\ve\ast V$ of the Riesz kernel $V(x)={\cal I}_{2s}(x)$ as
$$
V_\ve(x)= C_1\int \frac{\rho_\ve(y)}{|x-y|^{N-2s}}\,dy\,,
$$
where $C_1$ is the constant $C(N,s)$ of \eqref{riesz1}. Note the scaling
$$
V_{\ve\mu}(x) = V_{\ve}(x/\mu)/\mu^{N-2s}
$$
that allows to reduce the calculations to the case $\ve=1$. Clearly $V_\ve\to V$ in $L^1_{loc}(\ren)$. We want to examine further the approximation of $V$ by $V_\ve$.

\noindent $\bullet$ We claim that
$$
V_\ve(x)\ge  V(x) \qquad  \forall  |x|\ge \ve.
$$
Proof: It is enough to prove it for $\ve=1$. Using spherical coordinates we write $y=r\sigma$, $|\sigma=1|$,  $\int_{|x|\le 1}dx= \omega_N$, and then
$$
V_1(x)= C_1 \int_0^R \rho(r)\,dr \int_{|y|=r} \frac1{|x-y|^{N-2s}}\,r^{n-1}d\sigma=
N\omega_N C_1 \int_0^R \rho(r)r^{N-1} \Phi(r) dr
$$
where $\Phi(r):=\fint_{|y|=r} |x-y|^{-N+2s}\,d\sigma$.
So we are reduced to study the properties of the last integral. This is where we use the fact that for $z\ne 0$
$$
\Delta V(z)=(N-2s)(2-2s)V(z)/|z|^2\ge 0,
$$
so the integrand function is subharmonic in $B_r(x)$ if $r<|x|$. We conclude that
$$
\Phi(r):=\fint_{|y|=r}\frac1{|x-y|^{N-2s}}\,d\sigma\ge \Phi(0)=\frac1{|x|^{N-2s}}\,,
$$
and the result is true. In fact, Evans \cite{Evans}, page 26, proves that $$
\Phi'(r)=(r/N)\fint_{B_r(0)} \Delta V(x-y)\,dy=\frac{K}{r^{N-1}} \int_{B_r(0)} \Delta V(x-y)\,dy
$$
 hence $\Phi'(r)>0$ for all $r<|x|$.

\medskip

\noindent $\bullet$ We now estimate the error for large $|x|\ge 2$. In that situation we have
$$
\Phi'(r)\le c_2 r |\Delta V(x)|, \qquad \Phi(r)-\Phi(0)\le c_3r^2V(x)|x|^{-2}\,,
$$
so that
$$
V_1(x)-V(x)\le c_4 \left(\int_0^1 \rho(r)r^{N+1}dr\right) V(x)|x|^{-2}= c_5V(x)|x|^{-2}\,.
$$

Using now the rescaling rule, we get interesting results  for $V_\ve$.

\begin{prop} \label{prop.convve} There is a constant $K>0$ depending on $N$ and $s$ such that
\begin{equation}
V(x)\le V_\ve(x)\le V(x)\left(1+ K\frac{\ve^2}{|x|^2}\right)
\end{equation}
holds for all $\ve>0$ and all $|x|\ge 2\ve$. There is another constant $C$ such that
\begin{equation}
V_\ve(x)\le C V(x)
\end{equation}
for all $\ve>0$ and all $x$. Finally, $V_\ve-V\to 0$ in $L^1(\ren)$.
\end{prop}
For the last part, note that for $x$ away from zero $|V_\ve (x)-V(x)|\le C\ve^2|x|^{-N-2+2s}$
which is integrable at infinity, hence the $L^1$ norm in that region is small if $\ve\to 0$. Near zero $|V_\ve (x)-V(x)|\le CV(x)$ which is also integrable, hence uniformly small in a small ball.

\medskip

\noindent {\bf 2}. The following construction supplies an alternative approach that may also be useful. We consider ${\hat V}_n(x)=\inf\{V(x), n\}$ for $n\ge 1$. Note that this is just a rescaling on ${\hat V}_1(x)$ of the form \ ${\hat V}_n(x)= n\,{\hat V}_1( n^{1/(N-2s)} x)$. Clearly, ${\hat V}_n(x)$ increases monotonically to $V(x)$  as $n\to\infty$, ${\hat V}_n-V\to 0$ in $L^1(\ren)$.

  \begin{prop}  Consider the function $\phi_n(x):= (-\Delta)^s {\hat V}_n(x)$. Then, $\int \phi_n(x)\,dx=\int (-\Delta)^sV\,dx=1$ and $\phi_n(x)>0$ everywhere in $\ren$. The sequence $\phi_n=(-\Delta)^s {\hat V}_n(x)$ is a positive approximation sequence to the Dirac delta.
 \end{prop}

 \noindent {\sl Proof.} (i) We concentrate on $\phi_1$. On one hand, ${\hat V}_1$ has a flat plateau in a ball $|x|\le R$ where it attains the maximum value $1$, and  by a simple calculation with the explicit formula for the operator we get $(-\Delta)^s {\hat V}_1(x)>0$, and $\phi_1$ is bounded and continuous for $|x|<R$. In order to calculate the value of $(-\Delta)^s {\hat V}_1(x)$ at the points $|x|>R$ where ${\hat V}_1(x)=V(x)$ we define
$$
W(x)=V(x)-{\hat V}_1(x)=(V(x)-1)_+.
$$
It is also easy to see that $(-\Delta)^s W<0$ at the points $|x|\ge R$ since $W$ attains there its minimum, $W(x)=0$. The integral defining  $(-\Delta)^s W$ is bounded for $|x|>R$ and there is an asymptotic estimate $|(-\Delta)^s W|\sim C|x|^{-(n+2s)}$ as $|x|\to\infty$. Finally, by using integration with respect to a suitable cutoff function we prove that $(-\Delta)^s W$ is integrable and $\int (-\Delta)^s W (x) \,dx=0$.

Hence, $\phi_1(x)=(-\Delta)^s {\hat V}_1(x)>0$ also at those points. We also get the estimate
$$
\phi_1(x) \sim C(|x|^{-(n+2s)}) \quad \text{as } \ |x|\to\infty.
$$
Besides, $\phi_1(|x|)$ is continuous, bounded and smooth but for a possible asymptote at $|x|=R$.
Also, $\int \phi_1(x)\,dx=\int (-\Delta)^sV\,dx=1$.

(ii) Consider now $\phi_n(x):= (-\Delta)^s {\hat V}_n(x)$. It follows that $\phi_n$ is positive everywhere in $\ren$, that $\phi_n(x)=n^{N/(N-2s)} \phi_1(x n^{1/(N-2s)})$. This means that $   \int \phi_n(x)\,dx=\int \phi_1(x)\,dx $ and thus $\phi_n$ is a suitable approximation of the Dirac delta. \qed

\section{Special construction of a fundamental solution}\label{sec.specex}

\noindent (i) In the existence part we use  approximation of the special initial data $M\delta(x)$ with a smooth mollifier sequence of the form $M\rho_\ve(x)=M\ve^{-N}\rho(x/\ve)$. Without loss of regularity we fix $M=1$ here (using mass scaling, formula \eqref{mass.red}). It is then clear that the scaling group, formula \eqref{scaling1}, transforms the solution $u_\ve$ with data $\rho_\ve$ into the solution $T_\lambda u_\ve$ with data $\rho_{\lambda \ve}$.
By uniqueness of solutions we conclude that
\begin{equation}\label{id.le}
T_\lambda u_\ve=u_{\lambda \ve}\,.
\end{equation}
This will be used below.

\medskip

(ii) {\sc The Potential Equation, Continuity}.  The corresponding solutions $U_\ve$ of the PE (with initial data $V_\ve$) are regular in space, since $(-\Delta)^{s}U_\ve=u_\ve\in C_x^{\alpha}$. Since
\begin{equation}\label{eq.pot2}
\partial_t U_\ve=-u_\ve^m
\end{equation}
they  are also $C^{1,\alpha}$ in time uniformly for $t\ge \tau>0$. They are also monotone in time for all $t>0$.

Before examining the limit, we will establish another important property, the continuity of the evolution orbit $U(\cdot,t)$ in $L^1(\ren)$.  In fact, for every $0<t'<t$ we have
$$
\int |U(t')-U(t)|\,dx\le \int_{t'}^t\int |U_t|\,dxdt=\int_{t'}^t\int |u^m|\,dxdt<\infty.
$$
The worst case happens for $t'=0$ and we have proved above that the last integral can be estimated as $O(t^{2s\alpha})$, that goes to zero as $t\to 0$. This calculation is valid for the approximate solutions $U_\ve$, that are classical, and also for their limits.

We have assumed $m\ge 1$. For $m<1$  using an argument like \eqref{m<1} we have
\begin{equation}
\int (U(x,0)-U(x,t))\phi\,dx\le Ct.
\end{equation}
So in this case we have a uniform  continuity control in $L^1_{loc}(\ren)$. The rest is similar.

\medskip

\noindent (iii) {\sc Convergence for the Potential Equation}.
We have proved the lemma that says that $U_\ve(x,0)\le CV(x)$, therefore, $CV(x)$ is a uniform upper bound for the initial data, and in view of \eqref{eq.pot2}  also for the solutions.

In view of the internal regularity, we can pass to the limit along a subsequence $\ve_n \to 0$ to obtain a limit $U^*(x,t)$ that is a solution of the potential equation, $0\le U^*\le CV$.

As for the initial trace, in view of the convergence $V_\ve(x)=U_\ve (x,0)$ to $V$ in $L^1(\ren)$ (Proposition \ref{prop.convve}), and the uniform continuity of the orbits $U_\ve(t)$ in $L^1(\ren)$, the limit has also continuous time-increments in $L^1$, and takes the initial data $V$ in the sense that
\begin{equation}
\|U^*(t)-V\|_1\to 0 \quad \text{as } \ t\to 0.
\end{equation}
 Since $\partial_t U^*\le 0$, $U^*(x,t) $ converges monotonically to $V(x)$ as  $t\to 0$ for fixed $x$.

It is clear that for positive times the following holds: $(-\Delta)^s U^*(x,t)$ must coincide  with  $\overline{u}$, the limit  of $u_\ve(x,t)=(-\Delta)^s U_\ve(x,t)$ along $\ve_n$ (in other words, the operator is closed).
Since $\overline{u}(t)\to \delta $ in ${\cal D}'$ this also implies that $U^*(t)\to V$ in ${\cal D}'$, i.e., $U^*(x,0)=V(x)$.

The fact that the limit $U^*(x,t)$ takes the initial data with uniform convergence  away from zero is a consequence of the monotone convergence theorem for continuous functions. Around zero the convergence holds in  the Marcinkiewicz space $M^p(\ren)$, $p=N/(N-2s)$, hence in $L^q(\ren)$, $1<q<p$.

Besides, the convergence as $\ve \to 0$ is uniform for $t\ge \tau>0$.

We have concluded that $U^*$ is a solution of the PE corresponding to a fundamental solution of the original equation \eqref{eq1}, and it takes the initial data $V(x)$ in a strong way. We now have to solve the question of uniqueness of this limit $U^*$.

\medskip

\section{A uniqueness result for fundamental solutions}
\label{sec.uniq}

 \begin{thm} \label{thm-uniq} The fundamental solution $U(x,t)$ of the potential equation  with initial data $V(x)$ is unique under the assumptions that  $U\ge 0$, $U$ is continuous in $Q=\RR^N\times (0,\infty)$, and has only an initial singularity at $(0,0)$ in the sense that $U$ is bounded for $t\ge \tau>0$, is continuous  as $t\to0$ zero for all $x\ne 0$, and $U(x,t)\to \infty$ as $(x,t)\to (0,0)$.
\end{thm}

\noindent {\sl Proof.} Take two solutions $U_1$ and $U_2$ of the potential equation with the same initial data $V(x)$. Consider $U_3=U_1 + \delta$ and consider the approximate solution $\widehat U_2(x,t)= U_{2}(x,t+\tau)$ for a small $\tau>0$.
Then $U_3$ is strictly larger than the continuous and bounded function $\widehat U_{2}$ at $t=0$. The first point of contact cannot happen near infinity because of the uniform decay of $U_1$ and $U_2$ that implies strictly separation for $U_3$ and $\widehat U_2$.

At any  contact point $(x_0,t_0)$, $t_0>0$, and since $U_3- \hat U_{2}\ge 0$ up to that time, we have
$$
(-\Delta)^s (U_3-\hat U_2) < 0, \qquad \text{i.\,e. } \ u_{1}<\hat u_2,
$$
using the representation formula for $(-\Delta)^s$. On the other hand, at this point
$$
0\ge \partial_t(U_{3}-\hat U_{2})=- u_1^m+ \hat u_2^m, \qquad u_1^m\ge \hat u_2^m.
$$
This allows to arrive at a contradiction. We conclude that $
U_1(x,t)+\delta\ge U_2(x,t+\tau)$ for all $\delta,\tau>0$, hence in the limit
$$
U_1(x,t) \ge U_2(x,t) \qquad \text{for all } \ (x,t)\in Q\,.
$$
This concludes the proof, after reversing the roles of $U_1$ and $U_2$.

We now know that we do not need to take subsequences in the limit as $t\to0$ of Section \ref{sec.specex}. Note that we could have restricted the existence time and consider
$Q_T=\RR^N\times (0,T)$ instead of $Q$, and then the conclusion is valid for $0\le t\le T$.


\section{Self-similarity of the Barenblatt solution}
\label{sec.selfsim}

We conclude here the construction of the Barenblatt solution and the proof of Theorem \ref{thm.Bs}
for $N>1$.

\noindent $\bullet$ Let us settle first the question of self-similarity.
Indeed, the initial function of the solution of the PE with data $U_0(x)=C|x|^{-(N-2s)}$ is invariant
under the scaling $T$ defined by $(TU_0)(x)= A\, U_0(Bx)$ if we choose  $A=B^{N-2s}$. If moreover we have uniqueness of such solution, then the scaling group introduced in \eqref{scalPE} will imply that
$$
TU(x,t)=B^{N-2s}\,U(Bx, \lambda t)  \qquad \text{with} \ \lambda =A^{m-1}B^{2sm}=B^{N(m-1)+2s}
$$
will be a solution with the same initial data. By the claimed uniqueness $TU\equiv U$, hence for all $B>0$, $x\in\ren$ and $t>0$ we have
$$
U(x,t)=\lambda ^{(N-2s)\beta}\,U(\lambda^{\beta} x, \lambda t),
$$
where we have used the value for $\beta$ given in the Introduction.
Fix now $t_1>0$ and let $\lambda =1/t_1$. We get
$$
U(x,t)=t_1^{(N-2s)\beta}\,U(x\,t_1^{-\beta},1)\,.
$$
Of course, $t_1$ is arbitrary and can be replaced by $t$.
Calling now $U(x,1)\equiv H(x)$ we get
$$
U(x,t)=t^{-(N-2s)\beta} \, H(x\,t^{-\beta})\,.
$$
Applying operator $(-\Delta)^s$ we get self-similarity for the corresponding fundamental solution of \eqref{eq1}:
$$
u^*(x,t)=t^{-(N-2s)\beta} \, (-\Delta)_x^s (H(x\,t^{-\beta}))= t^{-\alpha} F(xt^{-\beta})\,.
$$
This is a rather classical argument in the study of self-similarity.

\noindent $\bullet$ {\sl End of proof of the Theorem.} The proof that $u^*(x,t)$ takes the data in the sense of measures (and not only as distributions) follows easy from self-similarity. The fact that $F(r)$ is monotone nonincreasing in $r$ is a consequence of the  Aleksandrov reflection principle proved in Section \ref{sec.alek}. We know by qualitative theory that $F$ is always positive since $u$ is. Same argument for the H\"older continuity. \qed

\medskip


\begin{figure}[!h]\label{figure_Barriers}
	\centering
\includegraphics[width=0.7\textwidth]{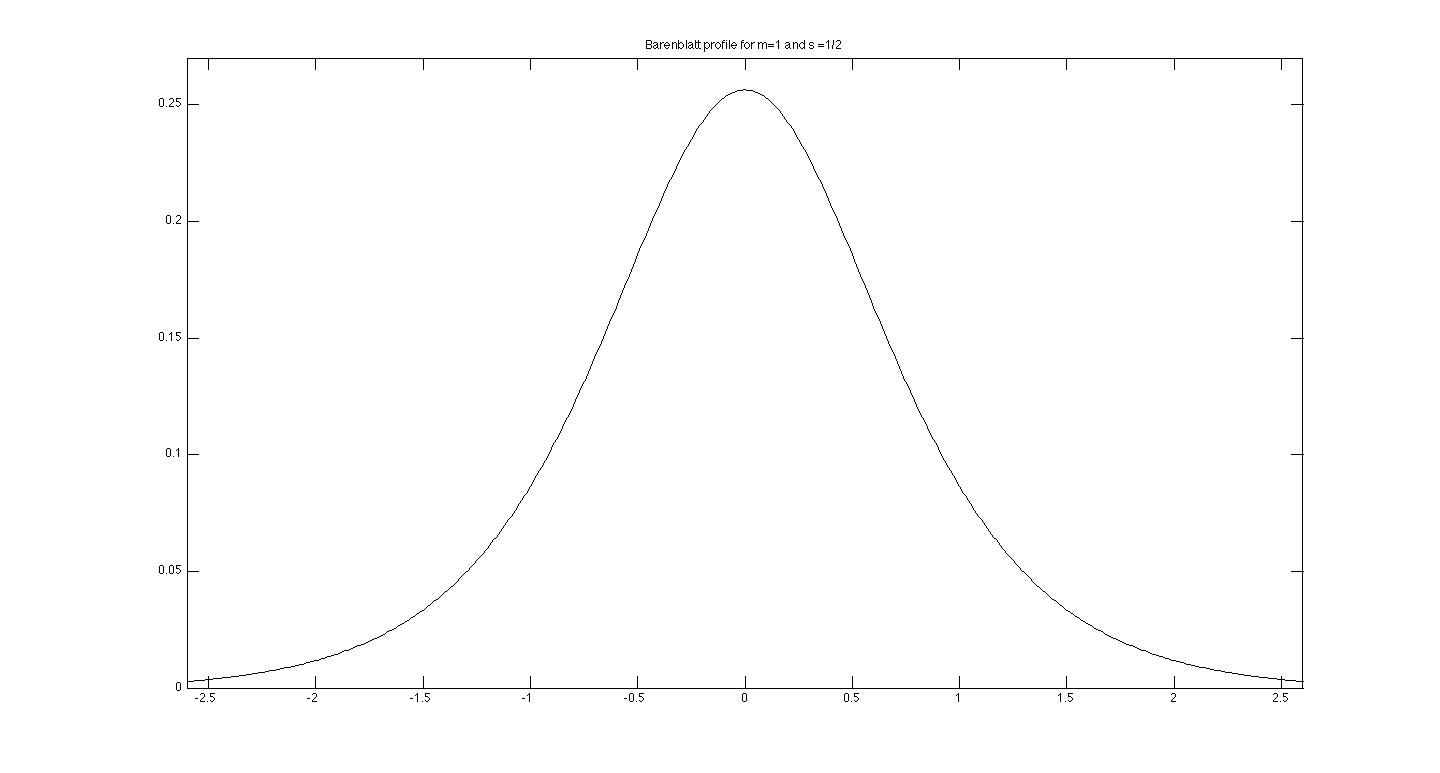}
\includegraphics[width=0.7\textwidth]{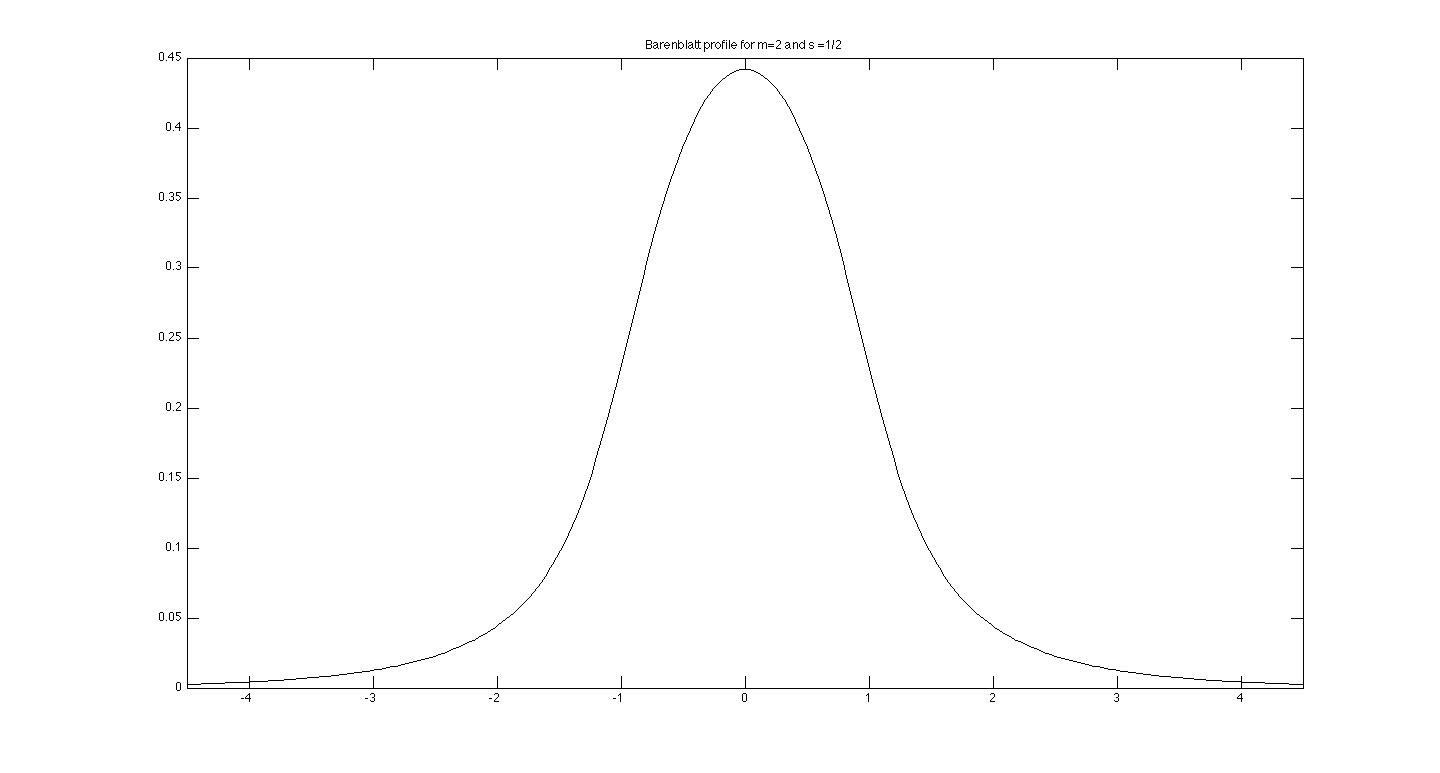}
\includegraphics[width=0.7\textwidth]{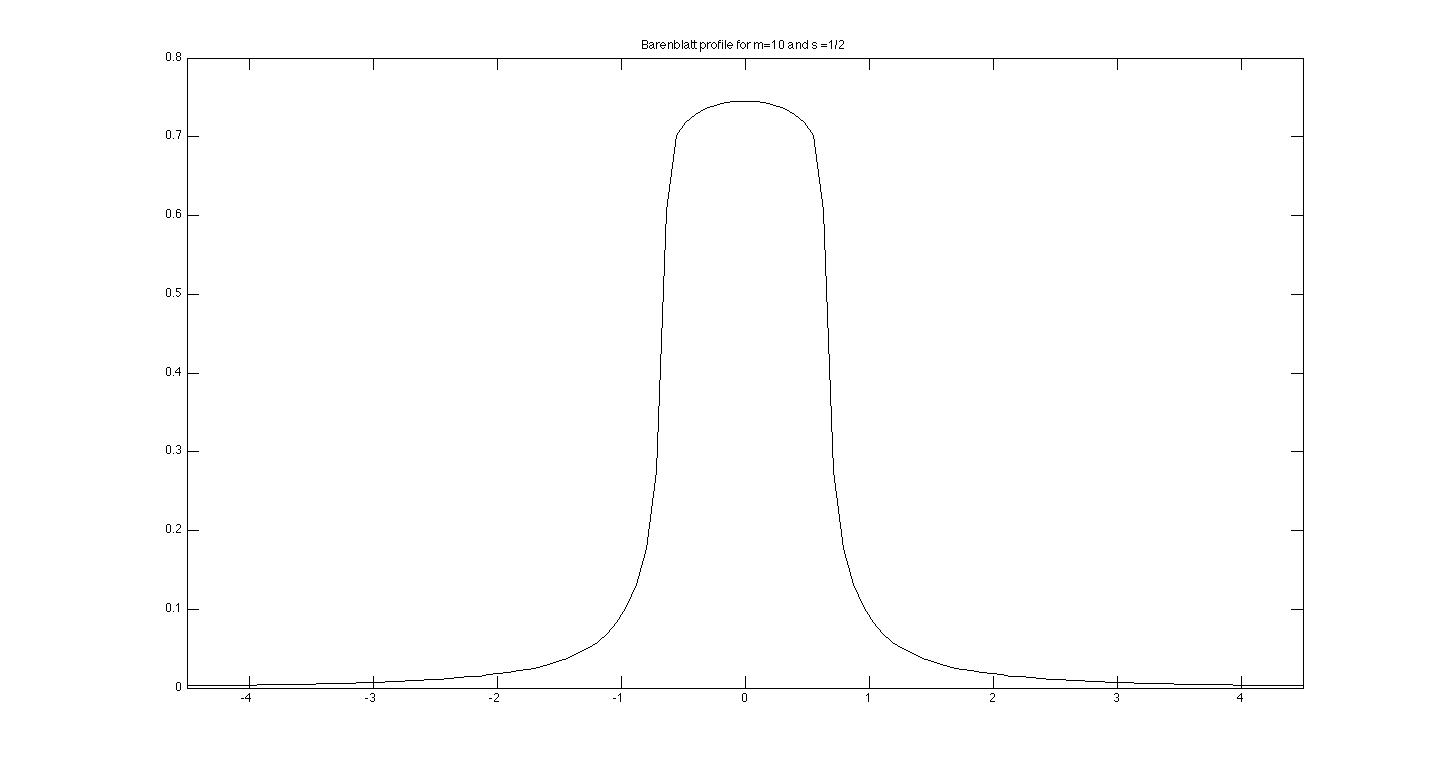}
\caption{Computed Barenblatt profiles for $m=1,2,10$ with $s=1/2$}
\end{figure}


\subsection{ Properties of the profile $F$}

\noindent {\bf Equation.} The self-similar profile $F$  satisfies an elliptic equation
\begin{equation}\label{sss.form}
(-\Delta )^{s}F^m =\alpha F +  \beta y\cdot \nabla F=\beta \nabla\cdot (yF)
\end{equation}
so that, putting $ s'=1-s$ we have
$$
\nabla (-\Delta )^{-s'} F^m=-\beta \, y\, F\,,
$$
which in radial coordinates gives
$$
L_{s'} F^m(r)=\beta \int_r^\infty rF(r)dr\,,
$$
where $L_{s'}$ the radial expression of operator $(-\Delta)^{-{s'}}$.

\medskip

\noindent {\bf Dependence on the mass.}  The scaling group acts on the profiles $F_M(r)$ for different masses $M>0$ and indeed we have
\begin{equation}
F_M(r)=\mu^{2s}F_1(\mu^{1-m} r), \quad M=\mu^{N(m-1)+2s}\,.
\end{equation}
which reduces all calculations to the case $M=1$. Since $N(m-1)+2s>0$ for $M>m_c$ we get $F_M(0)\to \infty$ as $M\to\infty$. For $m\ge1$ the same result happens for all $r>0$, $\lim_{M\to\infty} F_M(r)=\infty$. However, this last limit may be finite for $m<1$, see Section \ref{sec.vss} to understand when and why.

\medskip

\noindent {\bf Decay at infinity. First estimate.}
The precise behaviour of the fundamental profiles $F(y)=F_{m,N,s}(y)$ as $y\to\infty$ is a very important question in the qualitative theory. It is known in the linear case $m=1$, since $F$ is given by a linear kernel $K$ that decays like $|y|^{-(N+2s)}$, \cite{Blumenthal-Getoor}. The exact rate of decay for $m\ne 1$ is a nontrivial issue that we discuss next.

To begin with, the fact that $F$ is monotone as a function of $r$ and also integrable in $\ren$, implies that there is a constant $C=C(F)>0$ such that
  \begin{equation}
  F(r)  \le C\, r^{-N}\,.
  \end{equation}
We will use the same letter $C$ for different positive constants as long as their value is not important in the context. We will continue with this issue in Section \ref{sec.decay}.

\section{Peculiarities of one-dimensional flow}\label{sec.1d}

\noindent {\bf \ref{sec.1d}.1}. We add here the comments that are needed to close the case $N=1$, in the range $1/2<s<1$. In that case the kernel is a function that grows at infinity and the potential approach cannot be used in the direct way we have done before. There are at least 3 natural ways of addressing the difficulty.

One of them is to consider in a first stage the problem posed in a bounded domain $\Omega$ with nice boundary, say a ball $B_R$. Then the kernel is replaced by the Green function (with zero outside conditions) and this avoids considering the divergence of the kernel as $|x|\to\infty.$ Once the fundamental solution  for this problem is constructed, one passes to the limit as $R\to\infty$ and then proves uniqueness.

A second approach is to replace the fractional operator $(-\Delta)^{s}$  by a coercive operator like $L_\ve= (-\Delta)^{-s}+\ve I$, whose inverse has a nice kernel (the Bessel kernel $G_{2s}(x)$). If $U$ solves the potential equation as before and $L_\ve U=u$, then $u$ solves $u_t=L_\ve U_t=-L_\ve (u^m), $ i.\,e.,
\begin{equation}
u_t+(-\Delta)^su^m+\ve u^m=0.
\end{equation}
The plan is to construct a unique fundamental solution for this equation, and then pass to the limit and prove uniqueness.

The third option is to integrate in $x$, $v(x,t)=   \int_{-\infty}^x u(y,t)\,dy$. Then $v$ solves
\begin{equation}
v_t=L((v_x)^m), \qquad L= \partial_x(-\Delta)^{s-1}\,.
\end{equation}
THis is a kind of fractional $p$-Laplacian operator.

\medskip

\noindent {\bf \ref{sec.1d}.2}. Let us develop the second method. We take $\ve>0$ and solve the equation  $ u_t+(-\Delta)^su^m+\ve u^m=0$ for $x\in\re$ and $t>0$ with initial data $u_0\in L^1(\re)$ by just copying the method used in paper \cite{PQRV2} when $\ve=0$ (since the new term is dissipative, the needed estimates still hold). We obtain a solution with similar properties (except for the mass conservation). Again, for fixed $\ve>0$ the maps $S_\ve(t): u_0\mapsto u_\ve(t)$ are ordered contractions in $L^1(\re)$. If $u_0\ge 0$ so is the solution $u_\ve$ for all times and the family of solutions $\{u_\ve: \ve>0\}$ is monotone increasing as $\ve$ goes down to zero. In this way the standard solution $u$ of the problem  $ u_t+(-\Delta)^s(u^m)=0$, $u(x,0)=u_0$ is obtained in the limit $\ve\to 0$.

We now introduce the potential functions $U_\ve$ using the Bessel potential instead of the Riesz potential, so that $U_\ve(\cdot,t)=((-\Delta)^s +\ve I)^{-1}u_\ve(\cdot,t)$. The $U_\ve$ will satisfy, all of them, the  potential equation $U_t=-u^m$. We prove the existence and uniqueness of a solution with data $G_{2s}(x)=((-\Delta)^s +\ve I)^{-1}M\delta(x)$ much as above. This implies that the fundamental solution $u_{M,\ve}(x,t)$ of equation $ u_t+(-\Delta)^su^m+\ve u^m=0$ exists and is unique. We leave these details to the reader.

It follows easily that $u_{M,\ve}(x,t)$ is monotone increasing as $\ve$ goes down to zero and the limit is a fundamental solution of Equation \eqref{eq1} for $N=1$. We then prove that the limit is a minimal element among such possible fundamental solutions, and then conservation of mass implies that the fundamental solution must be unique. The result of Section \ref{sec.selfsim} follows then by easy changes in the argument.

\medskip

\noindent {\bf Note. } We have collected some properties of the Bessel potentials in an Appendix at the end of the paper for the reader's convenience. Note that these potentials have exponential decay as $|x|\to\infty$. Also (in $N=1$)  they are singular unbounded  at $x=0$ for $2s=\alpha<1$, but they are bounded for $2s\ge 1$.


\section{Asymptotic behaviour of general  solutions}\label{sec.ab}

We use the methods of the monograph  \cite{vazquez07} to get the following general theorem.
We recall that all our solutions are nonnegative.

\begin{thm}\label{thm.exlimit} Let $u_0=\mu \in {\cal M}_+(\ren)$, let $M=\mu(\ren)$ and let $u^* _M$ be the self-similar Barenblatt solution with mass $M$. Then  we have
\begin{equation}\label{conver.express}
\lim_{t\to\infty} t^{\alpha}\,|u(x,t)-u^*_M(x,t;M)|=0\,,
\end{equation}
and the convergence is uniform in $\ren$.
\end{thm}

\noindent {\sl Proof.} {\bf I.}   Firstly, we make the proof under the assumption that $d\mu=u_0(x)\,dx$ where $u_0$ is a bounded function with compact support. We divide this proof into several steps.

(i) We  may assume that there are $K>0$ big enough and $\tau>0$ such that
$u_0(x)\le u^*_K(x,\tau).$ The Maximum Principle implies then that
$$
u(x,t)\le u^*_K(x,t+\tau) \qquad \text{ for all } \ x\in\ren,  \ t>0..
$$
We now perform  the mass-preserving scaling transformation
$$
(T_\lambda u)(x,t)=\lambda^{\alpha} u(x\,\lambda^{\beta}, \lambda\, t)\,.
$$
It is easy to see that we have a similar estimate for the family $u_\lambda= T_\lambda u$:
$$
(T_\lambda u)(x,t)\le u^*_K(x,t+(\tau/\lambda)) \qquad \text{ for all } \ x\in\ren,  \ t>0\,.
$$
The family $(T_\lambda u)$ is bounded in $L^1\cap L^\infty\cap \,C^\alpha$ in any set of the form
$\Omega=B_R(0)\times (1/2,2)$ hence it converges to a solution $\widetilde u$ of the equation for $t>0$. The limit satisfies  the bound \ $\widetilde u(x,t)\le u^*_K(x,t)$  for all $ x\in\ren,  \ t>0\,$.

\noindent (ii) We claim that $\widetilde u$ takes the initial data $M\delta$ in the sense of measures. Indeed, away from zero the convergence of the limit $\widetilde u(x,t)$  to zero as $t\to 0$ is uniform, because it is bounded above by a tail of the self-similar Barenblatt solution $u^*_K(x,t)$. This leaves as ppssible initial trace a Dirac delta. It only remains to recall that the total mass of all the family $(T_\lambda u)$ is the same, $M$, at all times.

\noindent (iii) At this point we pass to the corresponding solutions $U_\lambda=(\Delta)^{-s}(T_\lambda u)$. They are solutions of the PE. The limit $U$ is also a solution of the PE with the already mentioned regularity and monotonicity properties.
The initial datum of $U(x,t)$ is $M$ times $V(x)$ in a weak sense, since
$$
\int (V(x)-U(x,t))\phi(x) \,dx =M(-\Delta)^{-s}\phi(0)-\int u(x,t)(-\Delta)^{-s}\phi(x)\,dx\,,
$$
which goes to zero since $\zeta(x)=(-\Delta)^{-s}\phi(x)$ is an acceptable test function for the convergence in the sense of measures of $u(\cdot,t)$ to $M\,\delta$.

On the other hand, the family $U(t)$ converges in $L^1(\ren)$ and monotonically to the initial trace that is $M\,V(x)$, as we have shown. The properties of the uniqueness theorem \eqref{thm-uniq} are met,
so that we conclude that $U$ is the fundamental solution $U_M^*(x,t)$.

\noindent (iv) Applying $(-\Delta)^{s}$ to $U=U^*$ we obtain the consequence $\widetilde u=u^*_M$. We thus have
$$
(T_\lambda u)(x,1)\to u^*_M(x,1) \qquad \text{uniformly in } \ x\in \ren.
$$
It is just a routine calculation to transform this expression into the convergence \eqref{conver.express}, see for instance  \cite{vazquez07}, chapter 18.

\medskip

\noindent {\bf II.} For a general initial function $u_0\in L^1(\ren)$ we have to work a bit more. First we fix an $\delta>0$ and truncate $u_0$ above and near infinity to fall into the previous case, namely
$$
0\le u_{0,\delta}(x) \le u_0(x), \qquad \|u_{0,\delta}\|_1=  M_\delta \ge M-\delta.
$$
In this way we have uniform convergence of the rescaled versions $u_\delta^{\lambda}(x,1)$ to $u^*_M(x,1; M_\delta)$ at $t=1$. Since the sequence  $\{ (T_\lambda u_\delta)(x,1): \delta>0\}$ is uniformly bounded and monotone in $\delta$ we get the convergence of the $\{(T_\lambda u_\delta)(x,1): \delta>0\}$ to $u^*_M(x,t;M)$ in $L^1$ and since they are bounded, they are uniformly $C^\alpha$, hence convergence is uniform on any fixed ball.

In order to examine the convergence of $(T_\lambda u_\delta)(x,1)$ outside a big  ball we first observe that all the functions are bounded in $L^1\cap L^\infty$, hence they are $C^\alpha$ with a uniform constant and exponent, which makes them a compact set of functions, hence we can extract a uniformly convergent subsequence. On the other hand, the $L^1$ norm of of $(T_\lambda u)(x,1)$ is small in the outer domain as a consequence of the following clever trick: fixing $\delta>0$ small, such $L^1$ norm  is bounded by  $\delta$ plus the norm of of $(T_\lambda u_\delta)(x,1)$ , and this one is controlled in that domain  by the tail of $u^*_K (x,1+\tau/\lambda)$, with $K=K(\delta)$ large. By interpolation between the $L^1$ and $C`^\alpha$ norms we get a small norm in $L^\infty$, uniformly in $\delta$ for $\ve$ small.

\medskip

(iii) For a measure as initial data, we displace the origin of times to $t=\tau>0$ and may then assume that $u_0$ is integrable and bounded.\qed

\medskip

\noindent  $\bullet $ We have a strong version of uniqueness of self-similar  solution, using the asymptotic behaviour:

\begin{cor}{\sl Let $u(x,t)$ be a nonnegative weak solution of Equation  \eqref{eq1}, and assume that $u$ is self-similar and the mass if finite, i.e., $\int u(x,t)\,dx<\infty$ for $t>0$. Then, $u$ is one of the Barenblatt solutions mentioned in the previous theorem, $u(x,t)=U_M(x,t)$, where $M=\int u(x,t)\,dx$, which is constant in time. }
\end{cor}

\section{Nonexistence  for small $m\le m_c$}

When we go  below the critical exponent $m_c$, fundamental solutions cease to exist. This was proved for the standard Fast Diffusion Equation ($s=1$, $0<m<(N-2)/N$) by Brezis and Friedman \cite{BF83}. Without entering in full details on the issue, we give here a simple proof of nonexistence in the fractional case, based on the existence of the scaling group \eqref{scaling1}. The main point to take into account is that when $m<m_c$ the similarity exponent $\beta$ of formula \eqref{scale.expo} is negative, i.\,e., the transformation involves space contraction instead of expansion. Therefore, if we consider a sequence of integrable functions approximating the Dirac delta of the form
$$
u_{0n}(x)=n^{N}u_{01}(nx), \quad \text{ with } \ \int u_{01}(x)\,dx=1,
 $$
 and say $u_{0n}\in C^\infty_c(\ren)$, then the sequence of solutions satisfy
$$
u_{n}(x,t)=n^{N}u_{1}(nx, n^{1/\beta}t)\,.
$$
Passing to the limit $n\to\infty$ for fixed $t>0$ the formula implies that $u_n(\cdot,t)\to \delta$ (thanks to the fact that $n^{1/\beta}\to 0$), and we conclude two things:  first, that we do not find in the limit the expected fundamental solution, and second, that the result can be interpreted as saying that, under this evolution equation, an initial Dirac delta does not spread with time and the ``physical solution'' is just $u_M^*(x,t)=M\delta(x)$.

In the critical case $m=m_c$ we have $1/\beta=0$ and the conclusion is the same.

We can try to take the limit of the fundamental solution that exists for $m>m_c$ and prove that the solution concentrates around the origin as $m\to m_c$. This is clear in the standard FDE since the solutions are explicit, see the very peculiar case $N=2$, $m\to 0$ in \cite{Vazsmooth06}, Lemma 8.3.

An analysis of the behaviour of the potentials is also illuminating. We leave it to the reader.

\section{Precise decay rates of the Barenblatt profiles}\label{sec.decay}

We resume the study of the decay rates of the Barenblatt solutions started at the end of Section \ref{sec.selfsim}.

 \noindent $\bullet$ {  \sc Case  $m>1$.} We may get the exact rate as follows: since $F^m\le Cr^{-mN}$ we obtain by use of the convolution formula that $L_{ s'}F^m\sim r^{-N+2s'}$, and this means that $\int_r^\infty r F(r)dr$ behaves like $r^{-N+2-2s}$ as $r\to\infty$. But this also implies that the ``mass'' in an annulus given by the integral $\int_r^{2r} r F(r)dr$ behaves like $r^{-N+2-2s}$ (with a different constant). Using the monotonicity of $F$ we arrive at the consequence that  $F$ has the  decay rate
 \begin{equation}
 F(r) \sim Cr^{-N-2s}\,,
 \end{equation}
 valid for all $m>1$, which equals the decay rate for $m=1$. This decay for $s<1$ is in stark contrast with the case $s=1$, where the profiles of the Barenblatt solutions are {\sl compactly supported} in the same $m$-interval, $m>1$.

\medskip

 \noindent $\bullet$ {\sc Case $m_c<m<1$}. A power-like bound from below is obtained  as follows: We may start from the homogeneity estimate \cite{BenCr} that says that $(1-m)tu_t\le u$. In terms of the self-similar profile,  this just means that $-(1-m)\beta(N\,F + rF'(r))\le F$, hence
$$
\frac{-rF'(r)}{F(r)}\le N+\frac{1}{(1-m)\beta}=\frac{2s}{1-m}\,.
$$
Integration of this inequality gives the following lower bound, valid for all $r\ge 1$, all $s\in (0,1)$ and all $m>m_c(s,N)$:
\begin{equation}\label{fde.rate}
F(r)\ge C\, r^{-2s/(1-m)}\,.
\end{equation}
Moreover, the function $ J(r):= F(r)\, r^{-2s/(1-m)}$ is monotone non-decreasing with $r$, so that it has a limit as $r\to\infty$. The lower bound \eqref{fde.rate} is a good starting point since it is the exact decay rate for the standard diffusion case $s=1$ when $m_c<m<1$ with compactly supported initial data, \cite{HP, ChVa02}.

We will  obtain a similar upper bound in the range $m_c<m<m_1$ as a consequence of the VSS construction in Section \ref{sec.vss}, where the asymptotic constant $c_\infty=\lim_{r\to\infty} J(r)$ is calculated. See more details in Section \ref{sec.vss}. Note in passing that $2s/(1-m)>N$ precisely for $m>m_c(s,N)$. Summing up, the decay rate \eqref{fde.rate} is optimal for $m< m_1=N/(N+2)$, and a very precise rate is obtained.

\medskip

 \noindent $\bullet$  However, the rate \eqref{fde.rate} is very far from a realistic estimate for $m$ close to 1, since $2s/(1-m)\to \infty$, which is  not, loosely speaking, an admissible decay in the situation of fractional Laplacian diffusion.

 We look for a bound from above as follows. Since $F(r)\le C\,r^{-N}$ we have $F^m\le  C\,r^{-Nm}$ and the fact that $Nm<N$  implies that $L_{ s'}F^m\le C\,r^{-Nm+2s'}$, hence
 $$
 \int_r^{\infty} rF(r)\,dr\le C\,r^{-Nm+2s'}.
 $$
  Using the monotonicity of $F$ this means that $F(r)\le C\,r^{-mN-2s}$. But $Nm+2s>N$ for $m>m_c$, so that we have a gain of upper bound exponent from $\gamma_0=N$ to $\gamma_1=mN+2s$. Proceeding iteratively as long as the whole argument we have used is justifed, we get a sequence of increasing exponents $\gamma_k$ given by $\gamma_{k+1}=\gamma_k m+2s$, that converge to the fixed point of the iteration formula, i.\,e., $2s/(1-m)$. This is all right for $m\le m_1$. However, for $1>m>m_1$ we have $2sm/(1-m)>N$ so that after a number of steps $F(r)^m\le r^{-m\gamma_k}$ with $\gamma_k m>N$. In that case, the next line of the previous argument changes to give \ $L_{ s'}F^m\sim r^{-N+2s'}$, and following the same line of argument $F(r)\sim r^{-N-2s}$ for all large $r$,which is the desired conclusion.

\medskip

 \noindent $\bullet$  In the limit case $m=N/(N+2s)$, the iteration from above is not inrerrupted and we get
 $F(r)\le C\,r^{-N-2s+\ve}$ for every $\ve>0$. Regarding the iteration from below, the starting rate is $\gamma=N+2s$, so that the exponent of our lower bound for $F^m$ is \ $mN/(N+2s)=N$. In this case the convolution formula for $L_{s'}$ produces a logarithmic correction for $L_{r'}F^m $ that we have to take into account. The estimate is then
\begin{equation}
F(r)\ge c r^{-N-2s}\log (r), \quad \mbox{for all large } \ r; \ m=m_1.
\end{equation}

We may sum up the results as follows.

\begin{thm} For every $m>m_1=N/(N+2s)$ we have the asymptotic estimate
\begin{equation}
\lim_{r\to\infty} F_M(r)\,r^{N+2s}=C_1\,M^{\sigma},
\end{equation}
where $M=\int F(x)\,dx$, $C_1=  C_1(m,N,s)>0$, and $\sigma=(m-m_1)(n+2s)\beta$. On the other hand, for $m_c<m<m_1$, there is a constant $C_\infty(m,N,s)$ such that
\begin{equation}
\lim_{r\to\infty} F_M(r)\,r^{2s/(1-m)}=C_\infty.
\end{equation}
The case $m=m_1$ is borderline and has a logarithmic correction.
\end{thm}

The analysis done before assumed $M=1$, for $M\ne 1$ just use scaling. For $m\in (m_c,m_1)$ the fact that $C_\infty$ does not depend on $M$ is commented and explained in whole detail in Section \ref{sec.vss}, where the last elements needed to complete this proof are given.


\begin{figure}[!h]\label{figure_Rates}
	\centering
\includegraphics[width=0.7\textwidth]{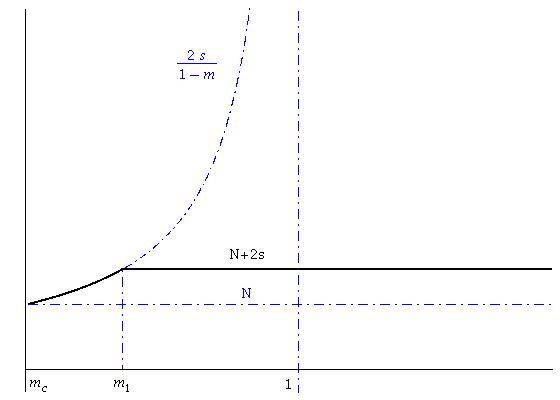}

\caption{Plot of decay rates of  the self-similar profiles. Here, $N=3$ and $s=0.8$.}
\end{figure}

\section{Very singular solutions in the fast diffusion range}\label{sec.vss}

\noindent {\bf \ref{sec.vss}.1. Existence.} A special solution that plays an important role in the theory of the standard Fast Diffusion Equation is the so-called  Very Singular Solution (VSS), cf.  \cite{ChVa02, Vazsmooth06}. More precisely, foer every $1>m>m_c$ there is a global in time function with the form of separation of variables that solves the FDE away from $x=0$ and has a standing singularity at $x=0$ for all times. It also starts with $V(x,0)=0$ for $x\ne 0$.

Using the same idea, we want to find  a similar solution for the fractional diffusion equation for $0<s<1$ in the good fast diffusion range $m_c<m<1$. We look for solutions of the form
$$
{\widetilde U}(x,t)= T(t)X(x)\,.
$$
Substitution into the equation leads to the value $T(t)=t^{1/(1-m)}$ for the time factor, while $Y=X^m$ has to satisfy
\begin{equation}\label{eq.Y}
(-\Delta)^s Y + \frac1{1-m} Y^p=0, \quad p=1/m.
\end{equation}
Now we try as solution the function $Y(x)= C^m |x|^{-\alpha} $. Under suitable conditions on $\alpha$, to be discussed below,
 we will have
\begin{equation}\label{lapl.power}
(-\Delta)^s Y = C^mk |x|^{-(\alpha+2s)}.
\end{equation}
for some negative constant $k=k(\alpha,N,s)$. It follows that Equation \eqref{eq.Y} for $Y$ can be satisfied only if $\alpha$ takes the value $\alpha(m,s)=2s/(p-1)=2sm/(1-m)$.
The constant $c>0$ is then determined by $
-C^mk(\alpha)=C/(1-m), $ so that
\begin{equation}\label{const.vss}
 C^{1-m}=(1-m)(-k(\alpha)).
\end{equation}
 Let us examine the conditions on $\alpha$.

(i) First, we need the condition $\alpha<N$ to make $Y$ integrable and formula \eqref{lapl.power} hold, which means $m<m_1:= N/(N+2s)$. Without this restriction we enter into a theory of Laplacians of non-integrable functions, which is acceptable anyway for $s=1$, but it does not work in the fractional case $s<1$, as we will see below.

(ii) We need to calculate $k$ to know when it is negative. This is done in \eqref{kalpha}. Note that  in the interval $m_c<m<m_1$ we have $N-2s<\alpha <N$ so that $N<\alpha+2s<N+2s$. We check there that $k(\alpha)>0$ for all $\alpha>N-2s$.

Finally, let us check the properties of the constructed function. We remark that ${\widetilde U}$ is not a weak solution in the standard sense since ${\widetilde U}(\cdot,t)$ is not even integrable in space around the origin. However, ${\widetilde U}^m$ is locally integrable, and we can propose a generalized  definition of weak solution that restricts  test functions to be supported away from the origin. We denote this situation with the name  {\sl generalized weak solution with an isolated singularity.}

\begin{thm} \label{thm.vss1} There exists a VSS for Equation \eqref{ell.1} in the range $m_c<m<m_1$ and it has the explicit form
\begin{equation}\label{vss}
{\widetilde U}(x,t)= C\,t^{1/(1-m)}\,|x|^{-2s/(1-m)}
\end{equation}
where $C=C(m,s,N)$  is  given by the explicit formula \eqref{const.vss}.
${\widetilde U}$ is a weak solution in a generalized sense,  more precisely, a weak solution with an isolated singularity.
\end{thm}

In order to compute the exact value of $C$ we must use \eqref{const.vss} together with \eqref{kalpha}. For further reference we callthe constant $C_{VSS}$.

\medskip

\noindent {\bf \ref{sec.vss}.2. VSS as limit.} We now establish the relation of the VSS with the Barenblatt solutions.

\begin{thm}\label{thm.vss2}
 The VSS \eqref{vss}  is the limit of the Barenblatt solutions $u^*_M(x,t)$ as the mass $M$ goes to infinity.
\end{thm}

\medskip

\noindent {\sl Proof.}
We  need to prove that this VSS  is the limit of the Barenblatt solutions as the mass goes to infinity. By the known comparison properties the sequence of Barenblatt solutions $u^*_M(x,t)$ is monotone increasing with $M>0$.

Next, we check that they are all bounded above by the VSS. This is done by comparison. But instead of comparing  $u^*_M$ with ${\widetilde U}$,  we compare $u(x,t)$ and ${\widetilde U}$, where $u$ is a solution with mass $M$ and smooth data $u_0$ supported in a small ball, so that $u(x,0)\le {\widetilde U}(x,1)$ for all $x\in\ren$. The comparison result holds and proves that $u(x,t)\le {\widetilde U}(x,t)$ for all $t>0$. Now we do the rescaling that is used in the proof of asymptotic behaviour and observe that ${\widetilde U}$ is invariant. Therefore,
$T_\lambda u(x,t)\le {\widetilde U}(x,t)$. In the limit $t\to\infty$ we use the result on asymptotic convergence to get $u^*_M(x,1)\le {\widetilde U}(x,1)$. The result is obviously true for all other times with a minimal modification of the argument.

This upper bound allows us to pass to the monotone limit  in the family $\{u^*_m(x,t): M>0\}$
and obtain a function $U^*(x,t)\le {\widetilde U}(x,t)$. This upper bound implies that $U^*$ is a weak solution of the Fractional FDE, and it is locally bounded away from the space origin $x=0$ for all times $t>0$. We conclude that $U^*$ is another possible very singular solution.

It also has the same invariance under $(T_\mu u)(x,t)=\mu^{2s/(1-m)}u(\mu x,t)$ which means that
the space form must be $c(t)|x|^{2s/(1-m)}$. But this implies the self-similar form we have chosen and this implies that $U^*={\widetilde U}$

\medskip

\noindent { {\bf \ref{sec.vss}.3. Precise asymptotic behaviour.} We claim that
the asymptotic constant
\begin{equation}
c_\infty(F_M)=\lim_{r\to\infty} F_M(r)r^{2s/(1-m)}
\end{equation}
for the Barenblatt solutions coincides with the constant $C$ for  VSS, let us call it $C_{VSS}$.

Indeed, by the ordering of the solutions we must have $c_\infty(F_M)\le C_{VSS}$. On the other hand, the scaling of the Barenblatt solutions immediately implies that $c_\infty(F_M)$ is independent of $M$. Therefore, for every $M>0$ we have
\begin{equation}
F_M(r)\le c_\infty\,r^{-2s/(1-m)}.
\end{equation}
If we had $c_  \infty(F) < C_{VSS}$, then the limit of the Barenblatt solutions as $M\to\infty$ would be a very singular solution of the form \eqref{vss} with a constant less than $C_{VSS}$, and this is not possible, as we have already seen.

\medskip

\noindent {\bf \ref{sec.vss}.4. Infinite limit in other ranges.} There are no fundamental solutions for $m\le m_c$ so that the ranges to discuss are $m\ge 1$ and $m_1\le m<1$. In the linear case $m=!$ the scaling of the Barenblatt solutions is
\begin{equation}
u^*_M(x,t)=M \,u(x, t)\,,
\end{equation}
so it is clear that $\lim_{M\to\infty} u_M^*(x,t)=\infty$ everywhere in $Q=\ren\times(0,\infty)$.
A similar result holds for $m> 1$,  in view of the scaling
\begin{equation}
u^*_M(x,t)=\mu^{2s}u(\mu^{1-m} x, t), \quad M=\mu^{N(m-1)+2s}\,.
\end{equation}
Since $u$ is monotone in $r=|x|$ and $m>1$, we have $u^*_M(x,t)\ge \mu^{2s}u( x, t)$ and the result follows letting $M\to\infty$.

In the remaining fast diffusion range, $m_1\le m<1$, the final answer is the same, an infinite limit, but the argument is not so simple. It is based on the minimum decay rate of the Barenblatt profiles, like $O(r^{-(N+2s)})$ (with logarithmic correction if $m=m_1$). In particular, such behaviour implies that the limit
$$
\lim_{ r\to\infty}F_1(r)t^{2s/(1-m)}=+\infty\,.
$$
The conclusion follows by scaling: For all large $M>0$ and putting $M=\mu^{N(m-1)+2s}$ we get for fixed $r>0$ as $M\to\infty$:
$$
\lim_{M\to\infty} F_M(r)=\lim_{\mu\to\infty} \mu^{2s}F_1(\mu^{1-m}r)=
r^{-2s/(1-m)} \lim_{z\to \infty} z^{2s}F_1(z)=\infty.
$$

\medskip

\noindent {\bf Remarks.} (1) This last divergence result means that no VSS in the sense of the beginning of this section can exist in this range. The result represents a marked difference with the case $s=1$, where a VSS solution exists for all $m\in (m_c,1)$, and the Barenblatt solutions are uniformly bounded by it.

   \medskip

\noindent  (2) The construction of a separated-variables solution with a fixed isolated singularity can de performed for $m<m_c$ but then we would get finite-time extinction profiles, in the line of \cite{Vazsmooth06}, Chapter \color{red} 5 \normalcolor. \normalcolor


\section{Construction of eternal solutions for $m=m_c$}\label{sec.eternal}

We have already seen that the Barenblatt solutions in the sense of fundamental solutions with self-similar form, do not exist for $m\le m_c:= (N-2s)/N$. One wonders if there is a special way to pass to the limit $m\to m_c$ to obtain a self-similar solution, even if it cannot be a fundamental solution.

\begin{thm} \label{thm.eter} There exists a family of eternal solutions of equation \eqref{eq1} with exponent $m=m_c$  and \ $0<s\le 1$. They have the form
\begin{equation}\label{vss.mc}
U(x,t)= e^{-Nct} F(|x|\,e^{-ct})
\end{equation}
where $F(r)$  is  a continuous, radial function with $F(0)=a>0$ and $F(r)\,r^N\to b>0$ as $r\to \infty$.
\end{thm}

\noindent {\sl Proof. \bf 1.} The process is known in the standard case $s=1$, but it will be convenient to explain it in some detail in this context, since it involves a subtle manipulation of the equations. We consider a different scaling of the flow that is good when passing to the limit $m\to m_c=(N-2)/N$ with $m>m_c$. It says
$$
u(x,t)= (1+ (ct/\beta))^{-n\beta} v(y, \tau), \quad y=x\,(1+ (ct/\beta))^{-\beta}, \ \tau=\beta \log(1+ c(t/\beta))\,.
$$
It applies for $t\ge -c\beta$, where $\beta$ is the self-similar exponent and $c>0$ is a free constant. The original equation\footnote{Note the slight change in the form of the equation due to the constant $m$, that is customary in studies of FDE. It is inessential here but we have inserted it for agreement with published information.}  $u_t=\Delta ( u^m/m)$  becomes
\begin{equation}\label{eq.renorm}
c\,v_\tau=\Delta (v^m/m)+ c\nabla \cdot (y\,v)=\nabla ( v\nabla (- p+ \frac{c}2 y^2)), \quad p=\frac1{1-m}v^{m-1}.
\end{equation}
Notice that no $\beta$ appears, which is important since $\beta\to\infty$ as $m\to m_c$. Equation \eqref{eq.renorm} has the stationary solution (profile)
$$
(\bar F_m(y))^{m-1}=(b+\frac{c(1-m)}2 y^2)\,,
$$
with $b>0$ arbitrary. This holds for $1>m>m_c$. In the limit $m\to m_c$ we have: $1/(1-m)\to N/2$, $\beta\to\infty$, $\tau(t)\to t$, $(1+ (ct/\beta))^{-\beta}\to e^{-ct}$,
$(1+ (ct/\beta))^{-N\beta}\to e^{-cNt}$, so that the profile in the limit is
$$
F_{m_c}(y)= (b+\frac{c}{N} y^2)^{-N/2},
$$
and the solution becomes (with $a=c/N$)
\begin{equation}
U_{m_c}(x,t)=e^{-cNt}(b + ax^2e^{-2ct})^{-N/2}=(ax^2 + b\,e^{2Nat})^{-N/2}\,,
\end{equation}
which is the (two-parameter family of) eternal solutions obtained in \cite{Vazsmooth06}, formula (5.36). We may still enlarge the family by displacing the origin of space coordinates. The scaling group $Tu(x,t) = \lambda^{N}\mu^{-N/2} u(\lambda x,\mu t)$, $\lambda, \mu>0$,  acts on these solutions to change the parameters $a$ and $b$ in the form: $T(a,b)=(a\mu, b\mu/\lambda^2)$, so the whole family of solutions comes from a single one.

Moreover, when $b\to 0$ we get the singular stationary case $U_\infty(x,t)= A|x|^{-N}$, $A=a^{-N/2}>0$, which is not a solution at $x=0$ since $\Delta U_\infty^m$ has a delta function there.  This situation reminds of the increasing limit of Barenblatt solutions to the VSS for $m>m_c$, in that case we had used the mass as paramter, here the parameter is rather $F_{m_c}(0)=b^{-N/2}$ and  $a>0$ is kept fixed. Note also that $U_\infty$ does not depend on time, and that the constant $A$ is arbitrary.

\medskip

\noindent {\bf 2}. We now apply these ideas  to the limit $m\to m_c(s)=(N-2s)/N$ in the fractional diffusion case $0<s<1$. This case is more difficult since the Barenblatt solutions are not explicit. We have to avoid the vanishing of the space asymptotic constant in the limit by careful time scaling of the above type that leads to the correct renormalized equation.

Indeed, the equations for $m>m_c$ have VSS solutions in the standard scaling of the form
$$
V(x)=C\,t^{1/(1-m)}x^{-2s/(1-m)}\,,
$$
 where $C_{VSS}$ is given by $C^{1-m}=-(1-m)k(\alpha)$ and
$$
 k(\alpha)=2^{2s}\frac{\Gamma((N-\alpha)/2)\,\Gamma((\alpha+2s)/2)}{\Gamma((N-\alpha-2s)/2)\,
\Gamma(\alpha/2)}
$$
with $\alpha=2ms/(1-m)$. When $m\to m_c$ we have  $\alpha\sim N-2s$, $(N-\alpha)/2\sim s$,
$(\alpha+2s)/2\sim N/2$, $\alpha/2\sim (N/2)-s$,  and finally
$$
\frac{N-\alpha-2s}2=\frac{N-2s}2-\frac{sm}{1-m}=\frac{N(1-m)-2s}{2(1-m)}=-\frac{1}{2(1-m)\beta}
$$
Therefore, $C^{1-m}\sim k_1/\beta$. Comparing the profile equations \eqref{sss.form} and
the one above, we see that we replace $\beta$ by $c$ so that the new constant $C$ is given by
$C^{1-m}=k_1/c$ which can be made 1 by suitable choice of $c$.

\medskip

\noindent {\bf  3. \sc Construction.} We now consider the family of self-similar solutions for $m\sim m_c$ with $F_m(0)=1$ (which means that you need to choose the mass of the solution appropriately). We also have
$$
F_m(r)\le \hat C_m r^{-2s/(1-m)}, \qquad \lim_{r\to\infty} F_m(r)r^{2s/(1-m)}=\hat C_m\to 1.
$$
Using the monotonicity of the profiles and $C^\alpha$ regularity of the family $F_m$ we may pass to the limit $F_m(r)\to F(r)$ so that $F(0)=1$, $F(r)\le  r^{-N},$ with $F\in C^\alpha$, and $F(r)\, r^{N}$ increasing. This is the eternal solution, after some extra work that is rather easy.

You may then construct a two parameter family of radial solutions, and then add  displacement in space. The two parameters are $F(0)$ controlled by constant mass scaling, and the constant at infinity that is controlled by the other part of the scaling group.\normalcolor \qed

\section{Aleksandrov's symmetry principle}\label{sec.alek}

The Aleksandrov-Serrin reflection method is a well-established tool to prove monotonicity of solutions of wide classes of (possibly nonlinear) elliptic and parabolic equations,  cf. \cite{Alek60, Ser71}. It has been quite useful in particular in the case of the PME, as documented in  \cite{vazquez07}. Here we will establish a version of the  principle valid in the presence of fractional operators, and then derive useful  monotonicity results.

\medskip

\noindent {\bf \ref{sec.alek}.1. \sc Elliptic setting.} We consider the problem
\begin{equation}\label{ell.1}
L u^m + u =f
\end{equation}
posed in $\ren$ with $L=(-\Delta)^s$, $m>0$, and $f\in L^1(\ren)$. We take a hyperplane $H$ that divides $\ren$ into two half-spaces $\Omega_1$ and $\Omega_2$ and consider the symmetry with respect to $\Pi$ that maps $\Omega_1$ onto $\Omega_2$. Then we have.

\begin{thm} Let $u$ be the unique solution of \eqref{ell.1} with data $f\in L^1(\ren)$, $f\ge 0$. Under the assumption that
\begin{equation}
f(x)\ge f(\Pi(x)) \qquad \text{in} \quad \Omega_1
\end{equation}
we have
\begin{equation}
u(x)\ge u(\Pi(x)) \qquad \text{in} \quad \Omega_1.
\end{equation}

\end{thm}

\noindent {\sl Proof.} (i) Due to the translation and rotation invariance of the equation we may assume that $H=\{x\in \ren: x_1=0\}$, and $\Omega_1=\{x\in \ren: x_1>0\}$. We write $x'=\Pi(x)$
where $\Pi $ is the symmetry with respect to  $x_1=0$. We take a function $f$  such that $f(x')\le f(x)$ when $x_1>0$. We assume at this stage that $f$ is bounded, continuous and integrable.

We solve the elliptic problem with data $f$ to get a solution $u$ and write $\widehat u(x)=u(x')$, that solves the problem with data $\widehat f(x)=f(x')$ in the whole space $\ren$.

\medskip

(ii) We now assune for a moment that $\widehat  u\le u$ in $\Omega_1$ and consider a point $x_0\in \Omega_1$ where $\widehat u$ touches $u$ from below. If such a point exists, then we have
$u(x')\le u(x)$ for all $x\in \Omega_1$ and we get that
$$
Lu^m(x_0)=C \int \frac{u^m(x_0)-u^m(y)}{|x_0-y|^{N+2s}}\,dy
$$
and
$$
L\widehat u^m(x_0)=C \int \frac{\widehat u^m(x_0)-C \widehat u^m(y)}{|x_0-y|^{N+2s}}\,dy=
C \int \frac{u^m(x_0)- u^m(y')}{|x_0-y|^{N+2s}}\,dy.
$$
Therefore, when calculating $L\widehat u^m(x_0)-Lu^m(x_0)$ the terms in $u^m(x_0)=\widehat u^m(x_0)$ cancel out and we get
$$
L\widehat u^m(x_0)-Lu^m(x_0)=C \int_{\Omega_1}\frac{u^m(y)-u^m(y')}{|x_0-y|^{N+2s}}\,dy+
C\int_{\Omega_2}\frac{u^m (y)-u^m(y')}{|x_0-y|^{N+2s}}\,dy\,.
$$
Interchange of $y$ and $y'$ in the last integral gives
$$
L\widehat u^m(x_0)-Lu^m(x_0)=\int_{\Omega_1}\frac{u^m(y)-u^m(y')}{|x_0-y|^{N+2s}}\,dy+
\int_{\Omega_1}\frac{u^m(y')-u^m(y)}{|x_0-y'|^{N+2s}}\,dy\,,
$$
hence,
$$
L (\widehat u^m-u^m)(x_0)=
\int_{\Omega_1}(u^m(y)-u^m(y')\left(\frac{1}{|x_0-y|^{N+2s}}-\frac{1}{|x_0-y'|^{N+2s}}\right)\,dy>0\,.
$$
where we have used that for $x_0,y\in \Omega_1$ we have $|x_0-y|\le |x_0-y'|$.
On the other hand, at this point we have $f(x_0)\ge f(x_0')$. We get a contradiction with the equation unless both solutions coincide everywhere. The solution must then equal its symmetric reflectio, which is impossible is $f$ is not symmetric.

\medskip

\noindent (iii)  In order to apply this argument we need to show that the situation where the solution and its reflection are ordered happens, which is not clear. The general argument that fits the actual situation is done by arguing on an approximation of the solution, as typical in the proofs of maximum principles for standard elliptic equations. We argue as follows:

- We take $f$ in, say, $L^\infty(\ren)\cap L^1(\ren)$ and with compact support. We use the theory of existence for data $\widehat f= f+C$, which is not difficult since it is like the one done in the papers \cite{PQRV1, PQRV2} after a vertical displacement of the solution. We thus get for bounded data $f_C=f+C$, $C>0$, a solution $u_C \ge C$ and $u_C- C\in L^1(\ren)$.

-  Under the stated hypotheses the existence and regularity theory says that  $u_C\in L^\infty(\ren)\cap L^1(\ren)$,  is continuous and that $u_C\ge  u_0= u$.
The regularity is obtained by a Morrey embedding theorem.

- A known strong maximum principle for fractional Laplacian operators  implies that these functions do not touch in $\ren$, i.\,e, $u_C >  u_0$. The argument is not very different from Step (ii) above, but simpler since there is no reflection.

- Let us now check the symmetrization argument on this family. We will write for more clarity $u_C = u(f+C)$, $u_0 = u(f)$. For large $C>0$ we have $u(f)\le u(f+C)$ and $\widehat {u(f)}\le u(f+C)$ in $\Omega_1$ ($\widehat {u(f)}$ indicates the symmetrical image of $u(f)$ as before).
We use the argument displayed before to conclude that there is no possible point of contact for $\widehat {u(f)}$ and $u(f+C)$ inside $\Omega_1$ and there was none at the boundary $H$ where
$\widehat{u(f)}=u(f)$.

- A this point we may lower the $C$ until we get a point of contact in one of the two comparisons. If $C>0$ there is no contact at $x_1=0$ and none as $|x|\to  \infty$ and we get a contradiction in the two comparisons.

 Therefore, the infimum of the $C$ is zero. Then, there is contact as $|x|\to  \infty$ and at the boundary $x_1=0$. But the plain comparison $\widehat {u(f)}\le u(f)$. The proof is finished under the stated assumptions on $f$.

 \medskip

 (iv) For general $f\in L^1(\ren)$ we prove the result by approximation with functions $f_n$ as above and then use the $L^1$-continuity of the map $f\mapsto u(f)$, \cite{PQRV2}.
 \qed

 \medskip

\noindent {\bf II. \sc Parabolic setting.} We consider the problem
\begin{equation}\label{par.1}
u_t+Lu^m=0,
\end{equation}
posed for  $x\in \ren$ and $t>0$ with initial data $u(x,0)=u_0(x)\in L^1(\ren)$, $u_0\ge 0$. We use the same notations and conventions for the hyperplane $H$, half-space $\Omega_1$ and symmetry $\Pi$.
Then we have

\begin{thm} Let $u$ the unique solution of \eqref{par.1} with initial data $u_0$. Under the assumption that
\begin{equation}
u_0(x)\ge u_0(\Pi(x)) \qquad \text{in} \quad \Omega_1
\end{equation}
we have for all $t>0$
\begin{equation}
u(x,t)\ge u(\Pi(x),t) \qquad \text{for} \quad x\in \Omega_1.
\end{equation}

\end{thm}

\noindent {\sl Proof.} This is an easy consequence of the fact that the nonlinear operator
$S: u\mapsto Lu^m$ can be defined as $m$-accretive in $l^1(\ren)$, and its resolvent $(I+\lambda S)^{-1}$ also and ordered contraction, cf. \cite{PQRV2}. We can then solve the parabolic problem by implicit discretization in time,
i.e., for times $t_0=0, t_1=h, t_2=2h, \cdots,$ we define the approximations
$$
\frac1h(u(t_{n}-u(t_{n-1}))+ L (u(t_n)^m)=0\,.
$$
We apply the conclusion  of the elliptic theorem iteratively at every step to get
$ u(x,t_n)\ge u(\Pi(x),t_n)$ for every $x\in \Omega_1$, $n\ge 1$. Finally, we pass to the limit
$h\to 0$ in the discretization step  (using the Crandall-Liggett Theorem as in \cite{PQRV2})
to get the solution $u(x,t)$ with the desired properties. See more details of this method in
\cite{CL, Cr86, vazquez07}. \qed

\medskip

\noindent {\bf III. \sc Application.}  Comparison after reflection around suitable hyperplanes
allows now to conclude that the solution with initial data supported in a small ball $B_{\ve}$
will be decreasing in the radial exterior direction in the annulus $\{x: |x|\ge 2\ve\}$. We also conclude that it is decreasing along a cone of directions centered along the radius with aperture that goes to the flat cone when the annulus is taken $\{x: |x|\ge R\}$ with $R/\ve\to\infty$. This implies that the solution with initial date a Dirac delta must be radially symmetric. All these techniques are well-known, as explained in  Chapter 9 of \cite{vazquez07}, where references to further literature are given.

\section*{Appendix 1. Fractional Laplacians and potentials}
\label{sec.app1}
\setcounter{section}{15}\setcounter{equation}{0}

According to  Stein,   \cite{Stein70}, chapter V, the definition of
$(-\Delta)^{\beta/2}$ is done by means of Fourier series
\begin{equation}
  ((-\Delta)^{\beta/2}f)^{\widehat{}}(x)=(2\pi|x|)^{\beta} \hat f (x)
\end{equation}
and can be used for positive and negative values of $\beta$. For
$\beta=-\alpha$ negative, with $0<\alpha<N$, we have the equivalence
with the Riesz potentials
\begin{equation}\label{riesz.formula}
      (-\Delta)^{-\alpha/2}f =I_{\alpha}(f)
      :=\frac1{\gamma(\alpha)}\int_{\mathbb{R}^N} \frac{f(y)}{|x-y|^{N-\alpha}}dy
\end{equation}
(acting on functions of the class $\cal S$ for instance) with
precise constant
$$
\gamma(\alpha) =\pi^{N/2}2^{\alpha}\Gamma(\alpha/2) /
\Gamma((N-\alpha)/2).
$$
 Note that $\gamma(\alpha)\to\infty$ as $\alpha\to N$, but $\gamma(\alpha)/(N-\alpha)$ converges to a
 nonzero constant, $\pi^{N/2}2^{N-1}\Gamma(N/2) $. Note that in this notation $\alpha$ stands for $2s$ in our previous notation.

\medskip

\noindent $\bullet$ The Fourier Transform of the function $f(x)=|x|^{-N+\alpha}$ is
 $\hat f(\xi)=\gamma(\alpha)(2\pi)^{-\alpha} |\xi|^{-\alpha}$.\\
The Fourier Transform of the function $f(x)=|x|^{-\alpha}$ is $\hat f(\xi)=\gamma(N-\alpha)(2\pi)^{\alpha-N} |\xi|^{\alpha-N}$.
We conclude that  its s-Laplacian, $(-\Delta)^s |x|^{-\alpha}$, is the power $k(\alpha)\,|x|^{-\alpha-2s},$ with a constant factor that equals
$$
k(\alpha)=\frac{\gamma(N-\alpha)}{\gamma(N-\alpha-2s)}=
\frac{\gamma(N-\alpha)}{\gamma(N-\alpha-2s)}
$$
so that
\begin{equation}\label{kalpha}
k(\alpha)=2^{2s}\frac{\Gamma((N-\alpha)/2)\,\Gamma((\alpha+2s)/2)}{\Gamma((N-\alpha-2s)/2)\,\Gamma(\alpha/2)}
\end{equation}
When $s=1$ we get
$$
k(\alpha)=\alpha(N-\alpha-2)
$$
which is positive whenever $\alpha\in (0,N-2)$. For $0<s<1$ it is not so clear.
In any case, in the range of interest, $0<N-2s<\alpha<n$, all the factors minus one are positive and the remaining one,
$\Gamma((N-\alpha-2s)/2)$ corresponds to an argument $((N-2s)-\alpha)/2$ which is negative but larger than $-s>-1$, hence $\Gamma((N-\alpha-2s)/2)<0$.

\section*{Appendix 2. Bessel kernels}

They are usually introduced via Fourier transform,
$$
\widehat G_\alpha(\xi)=(1+|\xi|^2)^{-\alpha/2}, \quad \xi\in\ren, \alpha>0.
$$
There is an expression for the kernel of the form
$$
G_\alpha(x) = c(\alpha)\int _0^\infty
e^{\pi |x|^2/t} \,t^{(\alpha-N)/2}\,\frac{dt}t, \qquad x\in \ren\,,
$$
 see \cite{Stein70}, Section V.3.1, so that $G_\alpha$ is a non-negative, radially decreasing function.  Moreover, $G_\alpha$ is integrable and  $c(\alpha)$ is a positive constant chosen so that  $\|G_\alpha\|_1=1$. The Bessel potential of order $\alpha>0$  of the density $\rho$ is defined by
$$
{\cal  B}_\alpha(x) =\int_{\ren} G_\alpha(x- y)   \rho(y)\, dy.
$$

We have the following estimates for  $0 < \alpha < N$:
$$
0 < G_\alpha(x) \le C\,|x|^{-(N-\alpha)}, \qquad \text{if} \quad 0<|x|<1,
$$
$$
 0 < G_\alpha(x) < e^{-|x|/2} \qquad \text{if} \quad |x|>1\,,
 $$
 where  $C = C(\alpha; N)$. More precisely, the kernel can  be represented by means of the McDonald function:
$$
G_\alpha(x) = c(\alpha, N) |x|^{(\alpha-N)/2} K_{(N-\alpha)/2} (|x|).
$$
The Macdonald function with index $\nu$, $K_\nu$, $\nu\in \re$, called also the
modified Bessel function of the second kind,  is given by the following
formula:
$$
K_\nu(r) = 2^{-1 -\nu} r^{\nu}\int_0^\infty  e^{-t} e^ {-r^2/4t}\,t^{-1-\nu}dt, \quad r > 0.
$$
The asymptotic behaviour of $K_\nu$ is as follows:
$$
K_\nu(r)\sim \frac{\Gamma(\nu)}2 (r/2)^{-\nu} \quad \text{as} \ r \to 0+, \ \nu>0
$$
while $K_0(r)\sim \log(1/r)$ as $r\to 0+$. At infinity we have
 $$
K_\nu(r)\sim \frac{\sqrt{\pi}}{\sqrt{2r}} e^{-r} \quad \text{as} \ r \to \infty,
$$
where the notation $g (r)\sim f (r)$ means that the ratio of $g$ and $f$ tends to 1. For $\nu < 0$ we have
$K_{-\nu}( r)= K_\nu( r),$ which determines the asymptotic behaviour for negative indices.
For properties of $K_\nu$ see \cite{Erdel}.

\section*{ Comments and extensions}

\noindent  $\bullet $ The linear case was known, and our results recover the information with a completely different machinery.

\noindent  $\bullet $ One wonders if there are any explicit or semi-explicit formulas for the family of Barenblatt solutions of this paper.

\noindent  $\bullet $ The limits $s\to 0$ and $s\to 1$ are worth examining.

\noindent  $\bullet $ Fundamental solutions can be constructed for the equation posed in a bounded domain. They do not play however such an important role in the theory. For instance, they shed no light on the asymptotic behaviour, even in the simplest case of zero Dirichlet boundary conditions.

\noindent  $\bullet $ The positivity of the fundamental solutions in all $\ren$, and more precisely, the power-like maximal decay rate as $|x|\to\infty$, are properties shared by all nonnegative solutions of the FPME and not only the solutions that take a special self-similar form. Quantitative versions of the lower bounds on the positivity of solutions are established in a paper with Bonforte \cite{BV12}. Higher regularity of positive solutions of equation \eqref{eq1}
is currently being investigated.

\noindent $\bullet $ There is the open problem of proving the uniqueness of solutions with measure data, not just Dirac deltas.

\noindent  $\bullet $ In the case of the standard PME, the range $m<m_c$ is now understood: no fundamental solutions exist, but the family of Barenblatt solutions for $m>m_c$ can be continued algebraically and they exhibit interesting attraction properties,   \cite{BBDGV, BDGVpnas}. Not much is known to our knowledge for the corresponding fractional equations.

\noindent  $\bullet $ New applications of these models in the applied sciences would be welcome as a motivation to develop further aspects of this theory.

\


\noindent {\large\bf Acknowledgments}

\noindent Work  partially supported by Spanish Projects MTM2008-06326-C0-01 and MTM2011-24696.
We thank Felix del Teso for the numerical computation of the Barenblatt profiles, which is part of his thesis work.

\vskip 1cm

\bibliographystyle{amsplain} 

\end{document}